\documentclass{article}

\usepackage{booktabs}
\usepackage{natbib}
\usepackage{multirow}
\usepackage{subcaption}
\usepackage{amssymb}
\usepackage{amsmath}
\usepackage{fullpage}
\usepackage{amsthm}
\usepackage{graphicx}
\usepackage{algorithm}
\usepackage{algpseudocode}
\usepackage{bm}
\usepackage{enumitem}

\usepackage{authblk}
\usepackage{adjustbox}

\usepackage{hyperref}
\usepackage{cleveref}

\hypersetup{
colorlinks=true,
linkcolor=red,
filecolor=magenta,      
urlcolor=cyan,
citecolor=blue,
}

\newtheorem{lemma}{Lemma}

\newtheorem{corollary}{Corollary}

\newtheorem{assumption}{Assumption}
\newtheorem{remark}{Remark}
\newtheorem{theorem}{Theorem}
\newtheorem{example}{Example}

\DeclareMathOperator*{\argmax}{arg\,max}

\DeclareMathOperator*{\app}{approx}

\newcommand{\p}{\mathcal P}

\newcommand{\1} { {(1)}}
\newcommand{\2} { {(2)}}

\newcommand{\E}{\mathbb E}

\newcommand{\h}{\mathcal H}

\newcommand{\dd}{{\mathfrak d}}
\newcommand{\pp}{{\mathfrak p}}

\newcommand{\dau}{\partial}
\newcommand{\what}[1]{\widehat{#1}}
\newcommand{\wtilde}[1]{\widetilde{#1}}

\newcommand{\upr}{{(r)}}
\newcommand{\nat}{\mathbb{N}_0}

\newcommand{\Unif}{\operatorname{Unif}}

\newcommand{\Mult}{\operatorname{Mult}}
\newcommand{\Mon}{\operatorname{Mon}}
\renewcommand{\Hat}{\operatorname{Hat}}

\newcommand{\R}{\mathbb{R}}

\newcommand{\ba}{\mathbf{a}}
\newcommand{\bb}{\mathbf{b}}

\newcommand{\mF}{\mathcal{F}}

\newcommand{\bzero}{{0}}

\renewcommand{\ba}{\mathbf{a}}
\newcommand{\bD}{\mathbf{D}}

\newcommand{\bp}{\mathbf{p}}
\newcommand{\bq}{\mathbf{q}}
\newcommand{\bu}{\mathbf{u}}
\newcommand{\bv}{\mathbf{v}}
\newcommand{\bw}{\mathbf{w}}
\newcommand{\bx}{\mathbf{x}}

\newcommand{\by}{\mathbf{y}}

\newcommand{\balpha}{\bm{\alpha}}

\newcommand{\bgamma}{\bm{\gamma}}

\newcommand{\bell}{{\bm{\ell}}}
\newcommand{\real}{{\mathbb R}}
\newcommand{\bsum}{{\beta_{\mathrm{sum}}}}
\newcommand{\rsum}{{\mathsf{r}_{\mathrm{sum}}}}
\newcommand{\eps}{\varepsilon}
\newcommand{\F}{\mathcal{F}}
\newcommand{\G}{\mathcal{G}}
\newcommand{\lt}{\left}
\newcommand{\rt}{\right}
\newcommand{\calP}{{\mathcal{P}}}
\newcommand{\calL}{{\mathcal{L}}}
\newcommand{\calN}{{\mathcal{N}}}

\newcommand{\calM}{{\mathcal{M}}}
\newcommand{\calB}{{\mathcal{B}}}
\newcommand{\calY}{{\mathcal{Y}}}
\newcommand{\calX}{{\mathcal{X}}}
\newcommand{\calZ}{{\mathcal{Z}}}
\newcommand{\bbN}{{\mathbb{N}}}
\newcommand{\beps}{{\boldsymbol{\epsilon}}}
\newcommand{\M}{\mathcal{M}}

\newcommand{\normal}{{\mathsf{N}}}

\algnewcommand\algorithmicinput{\textbf{INPUT:}}
\algnewcommand\INPUT{\item[\algorithmicinput]}
\algnewcommand\algorithmicoutput{\textbf{OUTPUT:}}
\algnewcommand\OUTPUT{\item[\algorithmicoutput]}

\title{A Likelihood Based Approach to Distribution Regression Using Conditional Deep Generative Models}

\date{\vspace{-3ex}}

\begin{document}
\author[1]{Shivam Kumar}
\author[2]{Yun Yang}
\author[3]{Lizhen Lin}
\affil[1]{Booth School of Business, University of Chicago}
\affil[2,3]{Department of Mathematics, University of Maryland, College Park}

\maketitle
\begin{abstract}

In this work, we explore the theoretical properties of conditional deep generative models under the statistical framework of distribution regression where the response variable lies in a high-dimensional ambient space but concentrates around a potentially lower-dimensional manifold. More specifically, we study the large-sample properties of a likelihood-based approach for estimating these models. Our results lead to convergence rate of a sieve maximum likelihood estimator (MLE) for estimating the conditional distribution (and its devolved counterpart) of the response given predictors in the Hellinger (Wasserstein) metric. 
Our rates depend solely on the intrinsic dimension and smoothness of the true conditional distribution. These findings provide an explanation of why conditional deep generative models can circumvent the curse of dimensionality from the perspective of statistical foundations and demonstrate that they can learn a broader class of nearly singular conditional distributions. Our analysis also emphasizes the importance of introducing a small noise perturbation to the data when they are supported sufficiently close to a manifold. Finally, in our numerical studies, we demonstrate the effective implementation of the proposed approach using both synthetic and real-world datasets, which also provide complementary validation to our theoretical findings.
\end{abstract}

\noindent\textbf{Keywords:} Distribution Regression; Conditional Deep Generative  Models; Intrinsic Manifold Structure; Sieve MLE;  Wasserstein Convergence.

\section{Introduction}\label{sec: Intro}
Conditional distribution estimation provides a principled framework for  characterizing the dependence relationship between a response variable $Y$ and  predictors $X$, with the primary goal of estimating the distribution of 
$Y$ conditional on $X$ through learning the (conditional) data-generating process. Conditional distribution estimation allows one to regress the entire distribution of 
$Y$ on $X$, which provides much richer information than the traditional mean regression  and plays a central role in various important areas ranging from causal inference \citep{causal, spirtes2010causal}, graphical models \citep{jordan1999graphical,koller2009graphical}, representation learning \citep{bengio2013representation}, dimension reduction \citep{carreira1997review, van2009dimensionality},  to model selection \citep{claeskens2008model, ando2010bayesian}.  Their applications span across diverse domains such as forecasting \citep{gneiting2014probabilistic}, biology \citep{krishnaswamy2014conditional}, energy \citep{jeon2012using}, astronomy \citep{zhao2021diagnostics}, and industrial engineering \citep{simar2015statistical}, among others.

There is a rich literature in statistics and machine learning on conditional distribution estimation including both frequentist and Bayesian methods  \citep{HallandYao2005, norets2017adaptive}. Traditional methods, both frequentist and Bayesian, suffer from the curse of dimensionality and often struggle to adapt to the intricacies of modern data types such as the ones with lower-dimensional manifold structures. 

Recent methodologies that leverage deep generative models have demonstrated significant advancements in complex data generation. Instead of explicitly modeling the data distribution, these approaches implicitly estimate it through learning the corresponding data sampling scheme. Commonly, these implicit distribution estimation approaches can be broadly categorized into three types. The first one is likelihood-based 
with notable examples including \citet{kingma2013VAE}, \citet{rezende2014stochastic}, \citet{burda2015importance}, and \citet{song2021maximum} . The second approach, based on adversarial learning, matches the empirical distribution of the data with a distribution estimator using an adversarial loss. 
Representative examples include \citet{GAN}, \citet{WGAN}, and \citet{sobolevGAN}, among others. The third approach, which is more recent, reduces the problem of distribution estimation to score estimation through certain time-discrete or continuous dynamical systems. 
The idea of score matching was first proposed in \citet{hyvarinen2005estimation} and \citet{vincent2011connection}. More recently, score-based diffusion models 
have achieved state-of-the-art performance in many applications \citep{sohl2015deep, nichol2021improved, song2020score, lipman2022flow, brehmer2020flows, han2022card}.

On the theoretical front, recent works such as \citet{huang2021wasserstein}, \citet{linsparse}, \citet{altekruger2023conditional}, \citet{stanczukdiffusion}, \citet{pidstrigach2022score} , and \citet{yunyangmanifold} demonstrate that distribution estimation based on deep generative models can adapt to the intrinsic geometry of the data, with convergence rates dependent on the intrinsic dimension of the data, thus potentially circumventing the curse of dimensionality.
Such advancement has naturally motivated us to  employ and investigate \emph{conditional deep generative model} for conditional distribution estimation.
Specifically, we explore and study the theoretical properties of a new likelihood-based approach to conditional sampling using deep generative models for data potentially residing on a low-dimensional manifold corrupted by full-dimensional noise.  
More concretely, we consider the following \emph{conditional distributional regression} problem:
\begin{equation}\label{eq: problem}
    Y|X = V|X + \eps,
\end{equation}
where $X$ serves as a predictor in $\real^\pp$, $V|X$ represents the (uncorrupted) underlying response supported on a manifold of dimension $\dd\le D$, $Y|X$ represents the observed response, and $\eps \sim \normal(0,\, \sigma_*^2 I_D)$ denotes the noise residing in the ambient space $\real^D$. Our deep generative model focuses on the conditional distribution $V|X$ by using a (conditional) generator of the form $G_*(Z,X)$, where $G_*$ is a function of a random seed $Z$ and the covariate information $X$. This approach is termed `conditional deep generative' because the conditional generator is modeled using deep neural networks (DNNs). 
Observe that, when $\dd < D$, the distribution of $G_*(Z, X)$ is supported on a lower-dimensional manifold, making it singular with respect to the Lebesgue measure in the $D$-dimensional ambient space. We study the statistical convergence rate of sieve MLEs in the conditional deep general model setup and investigate its dependence on the intrinsic dimension,  structure properties of the model as well as the noise level of the data.

\subsection{List of contributions}
We briefly summarize the main contributions made in this paper.
\begin{itemize}[leftmargin=*,itemsep=0pt]    
    \item To the best of our knowledge, our study is the first attempt to explore the likelihood-based approach for distributional regression using a conditional deep generative model, considering full-dimensional noise and the potential presence of  singular underlying support. 
    We provide a solid statistical foundation for the approach by proving the near-optimal convergence rates for this proposed estimator.

    \item We derive the convergence rates for the conditional density estimator of the corrupted data $Y$ with respect to the Hellinger distance and specialize the obtained rate for two popular deep neural network classes: the sparse and fully connected network classes. 
    Furthermore, we characterize the Wasserstein convergence rates for the induced intrinsic conditional distribution estimator 
    on the manifold (i.e.,~a deconvolution problem).  Both rates turn out to depend only on the intrinsic dimension and smoothness of the true conditional distribution.
    
    \item Our analysis in \Cref{cor: perturbation} suggests the need to inject a small amount of noise into the data when they are sufficiently close to the manifold. Intuitively, this observation validates the underlying structural challenges in related manifold estimation problems with noisy data, as outlined by \citet{Genovese}.
    
    \item We show that the class of learnable (conditional) distributions of our method is broad. 
    It encompasses not only the smooth distributions class, but also extends to the general (nearly) singular distributions with manifold structures, with minimal assumptions. 
\end{itemize}

\subsection{Other relevant literature}
The problem of nonparametric conditional density estimation has been extensively explored in statistical literature.
\citet{HallandYao2005}, \citet{bott2017nonparametric}, and \citet{bilodeau2023minimax} directly tackle this problem with smoothing and local polynomial-based methods. \citet{fan2004crossvalidation} and \citet{Efromovich2007regression} explore suitably transformed regression problems to address this challenge. Other notable approaches include the nearest neighbor method \citep{izbicki2020nnkcde, bhattacharya1990kernel}, basis function expansion \citep{sugiyama2010least, izbicki2016nonparametric}, tree-based boosting \citep{pospisil2018rfcde, gao2022lincde}, and Bayesian optimal transport flow \citet{chemseddine2024conditional} among others. 

In the context of conditional generation, we highlight recent work by \citet{huang2022sampling} and \citet{huang2021wasserstein}. In \citet{huang2022sampling}, GANs were employed to investigate conditional density estimation. While this work offers a consistent estimator, it lacks statistical rates or convergence analysis, and its focus is on a low-dimensional setup. 
In \citet{huang2021wasserstein}, conditional density estimation supported on a manifold using Wasserstein-GANs was examined. However, their setup does not account for smoothness across either covariates or responses, nor do they address how deep generative models specifically tackle the challenges of high-dimensionality.
Moreover, their assumption that the data lies exactly on the manifold can be restrictive. 
Our study shares some commonalities with the work of \citet{linsparse}, as both investigate sieve maximum likelihood estimators (MLEs). However, the fundamental problems addressed and the methodologies employed differ significantly, and our work involves technical challenges that span multiple scales. While \citet{linsparse} concentrates exclusively on unconditional distribution estimation, our theoretical analysis necessitates much more nuanced techniques due to the conditional nature of our setup. This shift is noteworthy because it demands a more refined analysis of entropy bounds, considering two potential sources of smoothness - across the regressor and the response variables. Furthermore, our setting accommodates the possibility of an infinite number of $x$ values, which gives rise to a dynamic manifold structure, further compounding the intricacy of the problem at hand.

\section{Conditional deep generative models for distribution regression}\label{sec: method}
We consider the following probabilistic conditional generative model, where for a given predictor value $x$, the response $Y$ is generated by
\vspace{-0.5em}
\begin{equation}\label{eq: model}
    Y  =  \, G_{*}(Z,x) + \eps, \quad x\in \calX \subset \mathbb R^\pp.
\end{equation}
Here, $G_*(\cdot,x) : \calZ \to \calM_x$ is the unknown generator function, $Z$  a latent variable with a known distribution $P_Z$ and support $\mathcal Z\subset \real^\dd$ independent of the predictor $X$. The existence of the generator $G_*$ directly follows from Noise Outsourcing \Cref{lemma: existence}. This lemma enables the transfer of randomness into the covariate and an orthogonal (independent) component through a generating function for any regression response. We denote $\calM:\, =\cup_{x\in \calX}\calM_x\subset \R^{D}$ as the support of the image of $G_*(\calZ, \calX)$  such as a (union of) $d$-dimensional manifold.  We model $G_*(\cdot, \cdot):\calZ \times \calX \subset \real^\dd\times \real^\pp\rightarrow \mathcal Y \subset \real^D$ using a deep neural network, leading to a \emph{conditional deep generative model} for \eqref{eq: model}. 

In the next section, we present a more general result in terms of the entropy bound (variance) for the true function class of $G_*$ and the approximability (bias) of the search class. We then proceed to a simplified understanding in the context of conditional deep generative models in subsequent sections.

\subsection{Convergence rates of  the Sieve MLE}\label{sec: theory}
In light of equation \eqref{eq: model}, it is evident that the distribution of $Y|X=x$ results from the convolution of two distinct distributions: the pushforward of $Z$ through $G_*$ with $X=x$, and $\eps$ following an independent $D$-dimensional normal distribution. The density corresponding to the true distribution $P_*(\cdot|X=x)$ can thus be expressed as:
\begin{equation}\nonumber
p_*(y|x) = \int \phi_{\sigma_*} (y - G_*(z,x)) \, dP_Z,
\end{equation}
where $\phi_{\sigma_*}$ is the density of $\mathsf{N}(\bzero,\sigma^2_* I_D)$.
We define the class of conditional distributions  $\calP$ as 
\begin{equation}\label{eq: P class}
\calP = \Big\{ P_{g,\sigma} : g(\cdot,x) \in \F, \sigma \in [\sigma_{\min}, \sigma_{\max}] \Big\},
\end{equation}
 where $P_{g,\sigma}$ represents the distribution with density $p_{g,\sigma} = \int \phi_\sigma (y - g(z,x))  dP_Z $. In this notation, \( P_* = P_{G_*, \sigma_*} \) and \( p_* = p_{G_*, \sigma_*} \). The elements of $\calP$ comprise two components: $g$ originating from the underlying function class $\F$, and $\sigma$, which characterizes the noise component. 
This class enables us to obtain separate estimates for $G_*$ and $\sigma_*$, furnishing us with both the canonical estimator for the distribution of $Y|X=x$ and enhancing our comprehension of the singular distribution of $G_*(Z, x)$, supported on a low-dimensional manifold.


Given a data set $\{(X_i, Y_i)\}_{i=1}^n$, the log-likelihood function is defined as $\ell_n(g,\sigma) = \tfrac{1}{n}\sum_{i=1}^n \log p_{g,\sigma}(Y_i|X_i)$. For a sequence $\eta_n \downarrow 0$ as $n \to \infty$, a \emph{sieve} maximum likelihood estimator (MLE) \citep{geman1982nonparametric} is any estimator $(\widehat{g}, \widehat{\sigma}) \in \F \times [\sigma_{\min}, \sigma_{\max}]$ that satisfies
\begin{equation}\label{eq: sieve MLE}
    \ell_n\lt(\widehat{g}, \widehat{\sigma}\rt) \ge \sup_{ \overset{\sigma \in [\sigma_{\min}, \sigma_{\max}]}{g \in \F} } \ell_n(g, \sigma) - \eta_n.
\end{equation}
Here $\what{g} \in \mathcal{F}$ and $\what{\sigma}\in [\sigma_{\min}, \sigma_{\max}]$ are the estimators, and $\eta_n$ represents the optimization error. The dependence of $\what{g}$ and $\what{\sigma}$ on $n$ illustrates the sieve's role in approximating the true distribution when optimization is performed over the class $\calP$. The estimated density $\what{p} = p_{\what{g}, \what{\sigma}}$ provides an estimator for $p_*(\cdot|\cdot)$, and $Q_{\what{g}}(\cdot | X=x)$ serves as the estimator for $Q_*(\cdot|X=x)$.

In this section, we formulate the main results, which provide convergence rates in the Hellinger distance for our sieve MLE estimator. The convergence rate was derived for any search functional class $\F$, with a brief emphasis on their entropy and approximation capabilities. 

\begin{assumption}[True distribution]\label{assume: basic}
Denote $\mu_X^\ast(x)$ as the distribution of $X$.
We   denote the true conditional densities as $p_*=\{p_*(\cdot|x), x\in \mathbb R^\pp\}$. It is natural to assume that 
the data is generated from $p_*$ from  model \eqref{eq: model} with some true generator $G_*$ and $\sigma_*$. 
We denote  $Q_*(\cdot|X=x)$ (or $Q_{G_*}$) as the distribution of $G_*(Z,x)$ for some distribution $P_Z$.
\end{assumption}

A function $g$ is said to have a composite structure  \citep{schmidtheiber, kohlerlanger} if it takes the form as
\begin{equation}\label{eq: composition}
    g = f_q \circ f_{q-1} \circ \cdots \circ f_1
\end{equation}
where $f_j : (a_j, b_j)^{d_j} \to (a_{j+1}, b_{j+1})^{d_{j+1}}$, $d_0 = \pp+\dd$ and $d_{q+1} = D$. Denote $f_j = (f_j^{(1)}, \ldots, f_j^{(d_{j+1})})$ as the components of $f_j$, let $t_j$ be the maximal number of variables on which each of the $f^{(i)}_j$ depends and let $f^{(i)}_j \in \h^{\beta_j}\lt((a_j, b_j)^{t_j}, K\rt)$ (see Section~\ref{sec-global} for the definition of the H\"{o}lder class $\h^\beta$). 
A composite structure is very general which includes smooth functions and additive structure as special cases. In addition,  in the next section, we show the class of conditional distributions $\left\{Q_{G_*}(\cdot|X=x): x\in\mathbb R^\pp, G_*\in \mathcal G\right\}$ induced by the composite structure is broad. 
\begin{assumption} [Composite structure]\label{assume: composition}
    Denote $\G = \G \lt(q, \bm{d},\bm{t}, \bm{\beta}, K \rt)$ as a  collection of functions of form \eqref{eq: composition}, where $\bm{d} = (d_0, \ldots, d_{q+1})$, $\bm{t} = (t_0, \ldots, t_{q+1})$, and $\bm{\beta} = (\beta_0, \ldots, \beta_{q+1})$. We regard $(q, \bm{d},\bm{t}, \bm{\beta}, K)$ as constants in our setup, and assume that the true generator $G_*(\cdot, x)$ as in \eqref{eq: model} belongs to $\G$, for all $x\in \calX$. Additionally, we assume $\||G_*|_\infty\|_\infty \le K$.    
        \begin{align*}
        &\wtilde{\beta}_j = \beta_j \prod_{l=j+1}^q \lt(\beta_l \wedge 1\rt), \quad j_* = \argmax_{j\in\{0, \ldots, q\}} \frac{t_j}{\wtilde{\beta}_j}, 
        \\
        &\beta_* = \wtilde{\beta}_{j_*}, \quad t_* = t_{j_*}.
    \end{align*}
The quantities $t_*$ and $\beta_*$ are called \textit{intrinsic dimension} and \textit{smoothness} of $G_*$ (or of $\G$).
\end{assumption}

\begin{remark}[Strength of the Composite structure]
    The expression $(a_j,b_j) \subset [-K,K]$ can be intuitively visualized by setting $a_j = -K$ and $b_j = K$. To illustrate the impact of intrinsic dimensionality and smoothness, consider a function $f:\mathbb{R}^d \to \mathbb{R}$ defined as $f(x) = f_1(x_1) + \ldots + f_d(x_d)$, where $x = (x_1, \ldots, x_d)$ and $f_j \in \mathcal{H}^\beta((-K,K), K)$ for $j = 1, \ldots, d$. While $f \in \mathcal{H}^\beta((-K,K)^d, K)$, its intrinsic dimension is $t_* = 1$ with intrinsic smoothness $\beta$. This mitigates the curse of dimensionality.
\end{remark}  
\begin{example}[One‐dimensional $\beta$-Hölder Generator]\label{ex:toy}
Let $U\sim\Unif(0,1)$ and define
\(
G(u)=u^{1/\beta}, u\in[0,1].
\)
Then $X=G(U)$ has density
\[
\frac{d}{dx} \mathbb{P}\bigl(U\le x^{\beta}\bigr)
=\beta\,x^{\beta-1},\quad x\in[0,1],
\]
which belongs to the $\beta$-Hölder class on $[0,1]$.  In our notation one checks
\(
t_*=1, \beta_*=\beta,
\)
and setting $\beta=1$ recovers $\Unif(0,1$) case and thus provides a fully explicit illustration of \Cref{assume: composition}.
\end{example}

\begin{assumption}\label{assume: manifold}
    Let $\calM_*$ be the closure of $G_*(\calZ, \calX)$. We assume that $\calM_*$ does not have an interior point, and reach$(\calM_*) = \mathsf{r}_\ast$ with $\mathsf{r}_\ast > 0$.
\end{assumption}

\Cref{assume: composition} permits low intrinsic dimensionality within the learnable function class. \Cref{assume: manifold} imposes the strong identifiability condition necessary for efficient estimation, as seen in manifold literature \citep{amarimanifold, yunyangmanifold}.


 Given two conditional densities $p_1(\cdot|x), p_2(\cdot|x)$ and $\mu_X^\ast$ denoting the density of $X$, we use integrated distances for a measure of evaluation. With a slight abuse of notation, we denote 
 $$
 d_1 (p_1, p_2) = \E_{X}\lt[ d_1 (p_1(\cdot|x), p_2(\cdot|x))\rt] \qquad \text{and}\qquad d_H (p_1, p_2) = \E_{X}\lt[ d_H (p_1(\cdot|x), p_2(\cdot|x)) \rt],
 $$ 
 where $d_1$ and $d_H$ represent the $L_1$ and the Hellinger distance as $d_1(p_1(\cdot|x), p_2(\cdot|x)) = \int \lt|p _1(y|x) - p_2(y|x) \rt| dy$ and $d_H (p_1, p_2) = ( \int\int [ \sqrt{p _1(y|x)} - \sqrt{p_2(y|x)} ]^2 \,  dy )^{1/2}$ respectively. 
Denote $\calN(\delta, \F, d)$ and $\calN_{[]}(\delta, \F, d)$ as covering and bracketing numbers of the function class $\F$ with respect to the (pseudo)-metric $d$. 

We first present \Cref{lemma: Entropy bound}, which  establishes the bracketing entropy of the functional class  $\mathcal P$ with respect to Hellinger
distance  in terms of the covering entropy of the search class $\F$.
This enables us to transfer the entropy control of the individual components $\F$ and $\sigma$ to the entire $\mathcal{P}$. 
\begin{lemma}\label{lemma: Entropy bound}
Let $\F$ be class of functions from $\calZ \times \calX$ to $\R^D$ such that $\||g|_\infty\|_\infty \le K$ for every $g \in \F$. Let $\calP = \lt\{ P_{g,\sigma}: g\in \F, \sigma \in [\sigma_{\min}, \sigma_{\max}] \rt\}$ with $\sigma_{\min} \le 1$. Then, there exist constants $c= c(\sigma_{\max}, K, D)$ and $C = C(\sigma_{\max}, K, D)$ and $\delta_* = \delta_*(D)$ such that for every $\delta \in (0, \delta_*]$,
    \begin{equation}\label{eq: lemma1}
        \adjustbox{max width = 0.9\linewidth}{$
        \log \calN_{[]}(\delta, \calP, d_H)  \le  \log \calN(c \sigma_{\min}^{D+3}\delta^4, \F, \||\cdot|_\infty \|_\infty) + \log \lt( \frac{C}{\sigma_{\min}^{D+2}\delta^4} \rt),
        $}
    \end{equation}
\end{lemma}
The proof of \Cref{lemma: Entropy bound} is provided in the \Cref{sec: proof lemma: Entropy bound}. 
\Cref{thm: main} presents the convergence rate of the sieve-MLE to the true distribution (see \Cref{sec: proof thm: main} for the proof). 
\begin{theorem}\label{thm: main}
Let $\F, \p,\sigma_{\min}$ and $\delta_*=\delta_*(D)$ be given as in \Cref{lemma: Entropy bound}, and $n \geq 1$. Suppose that 
$$
\log \calN \lt(\delta, \F, \|| \cdot |_\infty\|_\infty \rt) \le \xi \lt\{ A + 1 \vee \log \delta^{-1} \rt\}
$$
for every $\delta \in (0, \delta_*]$ and some $A,\xi >0$. Suppose that there exists a $G \in \F$ and some $\delta_{\app} \in (0, \delta_*]$ such that $\| | G - G_*|_\infty\|_\infty \le \delta_{\app}$. Furthermore, suppose that $s \ge 1$, $A \ge 1$, $\sigma_{min} \le 1$, $\delta_{\app} \le 1$ and $\sigma_* \in [\sigma_{\min}, \sigma_{\max}]$. Then \begin{equation}\label{eq: sieve MLE convergence sparse}
    P_* \lt( d_H(\widehat{p}, p_*) > \eps_n^* \rt) \le 5 e^{-C_1 n \eps_n^{*2}} + C_2n^{-1} 
\end{equation}
provided that $\eta_n \le n \eps^{*2}_n/6$ and $\eps^{*}_n \le \sqrt{2} \delta_{*}$, where
\begin{equation}\label{eq: rate}
    \eps_n^* = C_3 \lt( \sqrt{ \frac{\xi \lt\{ A + \log \lt( n/\sigma_{\min} \rt) \rt\}}{n}} \vee \frac{\delta_{\app}}{\sigma_*} \rt),
\end{equation}
$C_1$ is an absolute constant, $C_2 = C_2(D)$ and $C_3 = C_3(D, K, \sigma_{\max})$.
\end{theorem}
The outlined rate has two components: the statistical component, expressed as an upper bound to the metric entropy of $\F$, and the approximation component, denoted as $\delta_{\app}$. The statistical error is quantified by measuring the complexity of the class $\mathcal{P}$, as formulated in \Cref{lemma: Entropy bound}. The approximation error is assessed through the ability of the provided function class to approximate the true distribution.


\subsection{Neural network class}\label{se:NNC}
We model $G_*(\cdot, \cdot)$ using a deep neural network.
More specifically, 
we parameterize the true generator $G_*$  with a deep neural neural architecture $(L, {\bf r})$  of the form
\begin{equation}\label{eq: NN}
    \adjustbox{max width = 0.9\linewidth}{$
        f: \R^{r_0} \to \R^{r_{L+1}}, \quad z \mapsto f(z) = W_{L} \rho_{v_L} W_{L-1} \rho_{v_{L-1}} \ldots W_{1} \rho_{v_1} W_0 z,
    $}
\end{equation}
where $W_j \in \R^{r_{j+1} \times r_j}, v_j \in \R^{r_j}$, $\rho_{v_j}(\cdot) = \mathrm{ReLU}(\cdot - v_j)$ and ${\bf r} = (r_0, \ldots, r_{L+1}) \in \bbN^{L+2}$. The constant $L$ is the number of hidden layers and $r = (r_0, \ldots, r_{L+1})$ represents the number of nodes in each layer.

We define the \textbf{sparse} neural architecture class $\F_s(L,{\bf r}, s, B, K)$ as set of functions of form \eqref{eq: NN} satisfying 
\begin{equation*}
    \adjustbox{max width = \linewidth}{$
        \max_{0\le j \le L} |W_j|_\infty \vee |v_j|_\infty \le B, \quad \sum_{j=1}^L |W_j|_0 +  |v_j|_0 \le s, \quad \| |f|_\infty \|_\infty \le K,
    $}
\end{equation*}
with $r_0 = \dd + \pp$ and $r_{L+1} = D$, where $
|\cdot|_0$ and $|\cdot|_{\infty}$ stand for the $L^0$ and $L^\infty$ vector norms, and $ \| |f|_\infty \|_\infty = \sup_{x\in \real^{r_0} } \max_{i = 1, \ldots, D} |f_i(x)|$, $s$ is sparsity parameter and $K$ is functional bound. 

The \textbf{fully connected} neural architecture class $\F_c = \F_c\lt(L, {\bf r}, B, K\rt)$ is set of functions of form \eqref{eq: NN} satisfying
$$
\max_{0\le j \le L} |W_j|_\infty \vee |v_j|_\infty \le B, \qquad \||f|_\infty\|_\infty \le K.
$$
Both classes $\F_s$ and $\F_c$  for the deep generator will be considered in our analysis of  the resulting sieve maximum likelihood estimator. We denote the corresponding sieve-MLE  as $\widehat{p}_s$ and  $\widehat{p}_c$, respectively. When we use $r$ instead of ${\bf r}$, it refers to $r_1=\ldots = r_L = r$ along with $r_0 = \dd+\pp$ and $r_{L+1} = D$. 



We can simplify and visualize the result stated in \Cref{thm: main} in both cases: when the sieve-MLE is obtained with optimization performed over the class $\mathcal{F}_s$ and $\mathcal{F}_c$.
To fulfill the conditions stated in the \Cref{thm: main}, we need to establish entropy bounds for these function classes, $\mathcal{F}_s$ and $\mathcal{F}_c$, and gain insight into their approximation capabilities for the composite structure class described in \Cref{assume: composition}.

For the sparse neural architecture class $\mathcal{F}_s(L,r,s,K)$, the entropy, formally stated as Proposition 1 in \citet{ohn2019smooth}, is bounded as follows.
\begin{equation}\label{eq: entropy Hieber}
    \log \calN(\delta, \F_s, \||\cdot|_\infty\|_\infty ) \lesssim sL\, \{ \log(BLr) + \log \delta^{-1} \}.
\end{equation} 
From an entropy perspective, the fully connected neural architecture class $\F_c(L,r,B,K)$ can be viewed as $\F_s$ without any sparsity constraint, meaning $s \asymp r^2L$. Therefore, we have
\begin{equation}\label{eq: entropy Lin}
    \log \calN(\delta, \F_c, \||\cdot|_\infty\|_\infty ) \lesssim L^2 r^2 \{ \log(B L r) + \log \delta^{-1} \}.
\end{equation} 

The approximation properties of the sparse and fully connected network are provided in \Cref{lemma: DNN approximation}.1 and \Cref{lemma: DNN approximation}.2 of the \Cref{sec: approx old}, respectively.

Having established the essential components for $\mathcal{F}_c$ in \eqref{eq: entropy Lin} and Lemma \ref{lemma: DNN approximation}.2, and for $\mathcal{F}_s$ in \eqref{eq: entropy Hieber} and Lemma \ref{lemma: DNN approximation}.1, respectively, we can simplify Theorem \ref{thm: main} and state Corollary \ref{cor: both case}.

\begin{corollary}\label{cor: both case}
Suppose that Assumptions \ref{assume: basic} and \ref{assume: composition} hold, and $\sigma_* \in \left[ \sigma_{\min}, \sigma_{\max} \right]$ with $\sigma_{\min} \leq 1$ and $\sigma_{\max} < \infty$. Moreover, assume that the noise $\sigma_*$ decays at rate $\alpha$, i.e., $\sigma_* \asymp n^{-\alpha}$, and $\sigma_{\min} = n^{-\gamma}$ for some $\gamma \geq \alpha \geq 0$.
Then, for every $\delta_{\app} \in [0,1]$, the following holds:
\begin{enumerate}
    \item Let $\F_s = \F_s\lt(L, r, s, B, K \rt)$ with $\delta_* = \delta_*(D)$ be as given in \Cref{lemma: Entropy bound}, and $L \asymp \log \delta^{-1}_{\app} $, $r \asymp \delta^{-t_*/\beta_*}_{\app}$, $s \asymp \delta^{-t_*/\beta_*}_{\app} \log \delta^{-1}_{\app}$, $B \asymp \delta^{-1}_{\app}$. Then the sieve MLE $\what{p}_s$ satisfies \eqref{eq: sieve MLE convergence sparse} with $\eps_n^*$ as in \eqref{eq: rate} with $\xi = \delta^{-t_*/\beta_*}_{\app} \log^2 (\delta^{-1}_{\app})$ and $A = \log^2 (\delta^{-1}_{\app})$ provided that $\eta_n \le n \eps^{*2}_n/6$ and $\eps^{*}_n \le \sqrt{2} \delta_{*}$.
    
    \item Let $\F_c = \F_c\lt(L, r, B, K\rt)$ with $\delta_* = \delta_*(D)$ be as given in \Cref{lemma: Entropy bound}, and $L \asymp \log \delta^{-1}_{\app} $, $r \asymp \delta^{-t_*/2\beta_*}_{\app}$, $B \asymp \delta^{-1}_{\app}$. Then the sieve MLE $\what{p}_c$ satisfies \eqref{eq: sieve MLE convergence sparse} with  $\eps_n^*$ as in \eqref{eq: rate} with $\xi = \delta^{-t_*/\beta_*}_{\app} \log^2 (\delta^{-1}_{\app})$ and $A = \log^2 (\delta^{-1}_{\app})$ provided that $\eta_n \le n \eps^{*2}_n/6$ and $\eps^{*}_n \le \sqrt{2} \delta_{*}$.
\end{enumerate}
In particular, choosing
$
\delta_{\app} := \left( \sigma_*^2/n \right)^{\beta_*/(2\beta_* + t_*)} 
$
minimizes $\eps_n^* \asymp \sqrt{ \xi \lt\{ A + \log \lt( n/\sigma_{\min} \rt) \rt\}/n} \vee \delta_{\app}/\sigma_* $, and gives
\begin{equation}\label{eq: eps rate 2}
    \eps_n^* \asymp n^{-\frac{\beta_* - t_* \alpha}{2\beta_* + t_*}} \log^{2} (n).
\end{equation}
\end{corollary}

\begin{remark}
    The convergence rate in \eqref{eq: eps rate 2} illustrates the influence of intrinsic dimensionality, smoothness, and noise level on the estimation process. Note that $\alpha$ is upper bounded as $\eps_n^*\le\sqrt{2}\delta_*(D)$.  For large values of $\alpha$, estimation of $G_*$ is inherent difficult as the data is very close on the singular support. To address this, a small noise injection, as described in \Cref{cor: perturbation}, can smooth the estimation and ensure consistency.
\end{remark}
The proof of \Cref{cor: both case} is provided in \Cref{sec: proof cor: both case}. For the composite structural class $\G$, the effective smoothness is denoted by $\beta_*$, and the dimension is $t_*$. This effectively mitigates the curse of dimensionality. The convergence rate at \eqref{eq: eps rate 2} also recovers the optimal rate when $q=1$ and $\alpha =0$, and there is a small lag of polynomial factor $t_*\alpha/(2\beta_* + t_*)$ when $\alpha>0$ \citep{norets2017adaptive}. This lag arises due to the presence of full-dimensional noise in the response observation $Y$.
Note that when the noise is small, that is $\alpha$ is large, achieving a sharp estimation of $p_*$ requires an equally accurate estimate of $G_*$. This can be quite challenging. Our practically tractable approach attempts to address this without initially estimating the singular support.

\subsection{Wasserstein convergence of the intrinsic (conditional) distributions}
Using Wasserstein distance as a metric for distributions $Q_g$ is meaningful  due to their singularity in  ambient space: when $\dd<D$, the conditional distribution is singular with respect to the Lebesgue measure on $\real^D$.

The integrated  Wasserstein distance, for $r\ge1$, between $P_1(\cdot|X)$ and $P_2(\cdot|X)$ is defined as
\begin{equation*}
    \adjustbox{max width = \linewidth}{$
    W_r \lt( P_1, P_2 \rt) = \E_X\lt[\inf_{\beta\in \Gamma(P_1, P_2)} \lt( \E_{(U_1, U_2)\sim\beta}\big[ |U_1 - U_2|_r^r \big] \rt)^{1/r}\rt], 
    $}
\end{equation*}
where $\Gamma(P_1, P_2)$ is the set of all couplings between $P_1$ and $P_2$ that preserves the two marginals. The (dual) representation of this norm, $W_r(P_1, P_2) = \E_{X} \lt[ \sup_{\|f\|_{{\mathrm Lip}_r} \le 1} \lt\{ \, \E_{P_1} [f] -  \E_{P_2} [f] \,\rt\} \rt] $ \citep{villani2009optimal} with $\|\cdot\|_{{\mathrm Lip}_r}$ denoting the $r$-Lipschitz norm, is particularly useful in our proofs.
\begin{theorem}\label{thm: wasserstein general}
    Suppose that Assumption \ref{assume: manifold} holds. If $d_H(p_{g,\sigma}, p_*) \le \eps$ holds for some $\eps \in[0,1]$ and some $p_{g,\sigma} \in \calP$, then we have
    $$
    W_1(Q_g, Q_*) \le C \lt(\eps + \sigma_* \sqrt{\log \eps^{-1}} \rt),
    $$
    where $C= C(D, K, \mathsf{r}_\ast)$ depends only on $(D,K,\mathsf{r}_\ast)$.
\end{theorem}
The proof of \Cref{thm: wasserstein general} is provided in \Cref{sec: proof thm: wasserstein general}. \Cref{thm: wasserstein general} guarantees that $W_1\big(\what{Q}_{\what{g}}, Q_*\big) \lesssim_{\log} d_H(\what{p}, p_*) + \sigma_*$, where \( \lesssim_{\log} \) represents less than or equal up to a logarithmic factor of $n$. Following from \Cref{cor: both case}, the Wasserstein convergence rate, $n^{-(\beta_* - t_* \alpha)/(2\beta_* + t_*)} \log^{2} (n) \vee \sigma_* \log^{1/2}(n)$, comprises two components: the convergence rate in the Hellinger distance and the standard deviation of the true noise sequence. It is noteworthy that the first expression is influenced by the variance of noise by the factor $\alpha$. When $\alpha$ is very small, indicating that the data $Y_j$ lies very close to the manifold, the second expression $n^{-\alpha}$ in the overall rate dominates. Intuitively, this phenomenon arises from the underlying structural challenges in related manifold estimation problems with noisy data, as discussed by \citet{Genovese}. To address this issue, we propose a data perturbation strategy by transforming the data $\{(Y_j, X_j)\}_{j=1}^n$ into $\{(\wtilde{Y}_j, X_j)\}_{j=1}^n$, where $\wtilde{Y}_j = Y_j + \beps_j$ and $\beps_j \sim \normal\lt(0_D, n^{-\beta_*/(\beta_* + t_*)} \,  I_D \rt)$. The resulting estimation error bound is summarized below, whose proof is provided in \Cref{sec: proof cor: perturbation}.
\begin{corollary}\label{cor: perturbation}
    Suppose that Assumption \ref{assume: basic}, \ref{assume: composition}, and \ref{assume: manifold} hold, and $\sigma_* \in [\sigma_{\min}, \sigma_{\max}]$ with $\sigma_* = n^{-\alpha}$ and $\sigma_{\min} = n^{-\gamma}$ for some $0 \le \alpha \le \gamma$.  Then for each of the network architecture classes (sparse and fully connected) with the network parameters specified in \Cref{cor: both case}, the sieve MLE $\what{p}_{per}$ and $\what{Q}_{per}$ based on the perturbed data $\{(\wtilde{Y}_j, X_j)\}_{j=1}^n$ satisfies
\begin{equation*}
\adjustbox{max width=\linewidth}{$
    P_*\left[ W_1\left( \what{Q}_{per}, Q_*\right) \ge \left( \eps^*_n + \sigma_* \sqrt{\log \left((\eps^*_n)^{-1}\right)} \right) \right] \lesssim 5 e^{-C_1 n{\eps^*_n}^2} + \frac{C_2}{n}
$}
\end{equation*}
    where $\eps^*_n$ can be chosen such that
    \begin{equation}\label{eq: manifold rate}
        \adjustbox{max width=\linewidth}{$
            \eps^*_n + \sigma_* \sqrt{\log ((\eps^*_n)^{-1})}\, \asymp 
        \begin{cases}
        n^{-\frac{\beta_* - t_*\alpha}{2\beta_* + t_*}}\, \log^{2} (n) , \qquad &\text{if}\; \alpha < \beta_*/\{2(\beta_* + t_*)\} ,    
        \\
        n^{-\frac{\beta_*}{2(\beta_* + t_*)}}\, \log^{2} (n) , \qquad &\text{otherwise}.
        \end{cases}
        $}
\end{equation}
\end{corollary}

\subsection{Characterization of the learnable distribution class}
Section \ref{se:NNC} focuses on the true generator $G_*$ within the class of functions with composite structures. 
In this subsection, we show that such a conditional  distribution class achieved by the push-forward map
$G_*$ is  broad and includes many existing distribution classes for $Q_*$ as special cases. 

\subsubsection{Smooth  conditional density}
\label{sec-global}

For $\beta>0$, let $\h^\beta(D, M)$ be the class of all $\beta$-H\"older functions $f:D \subset \real^{\mathfrak d} \to \real$ with $\beta$-H\"older norm bounded by $M>0$. Let $\h^\beta(D) = \cup_{M>0} \h^\beta(D, M)$. See \Cref{sec: notations} for their formal definitions.
\begin{lemma}\label{lemma: Caffareli}
    Suppose that (i) $\calZ \times \calX $ and $\calY$ are uniformly convex and (ii) $p_Z \in \h^{\beta_Z}(\calZ)$, $\mu^*_X \in \h^{\beta_X}(\calX)$ and $q_* \in \h^{\beta_Q}(\calY)$ for some $\beta_Z, \beta_X, \beta_Q >0$ and are bounded above and below. 
    Then, there exists a map $g(\cdot,\cdot): \calZ \times \calX \to \calY$ such that $Q_*(\cdot|\cdot) = Q_g$ and $g \in \h^{\beta_{\min}+1}(\calZ\times \calX) $,
    where $\beta_{\min} = \min\{\beta_Z, \beta_X, \beta_Q \}$.
\end{lemma}
\Cref{lemma: Caffareli} establishes that the learnable distribution class includes Hölder-smooth functions with smoothness parameter $\beta_{\min}$ and intrinsic dimension $\dd$. As a result, following \Cref{cor: both case}, the convergence rate for density estimation is given by \(\eps^*_n \asymp n^{- (\beta_{\min} + 1 - \dd\alpha)/(2\beta_{\min} + 2 + \dd)}\).
A push-forward map is a transport map between two distributions. 
The well-established regularity theory of transport map in optimal transport is directly applicable here [see \citet{villani2009optimal} and \citet{villani2021topics}]. 
 The proof of \Cref{lemma: Caffareli} is based on Theorem 12.50 of \cite{villani2009optimal} and \citet{caffarelli1996boundary}, which establishes the regularity of this transport map  and its existence follows from  \citet{brenier1991polar}.
 When \( p_Z \) is selected as a well-behaved parametric distribution, the regularity of the transport map is determined by the smoothness of both \( \mu_X^* \) and \( Q_* \). For a more detailed discussion on this, please refer to \Cref{sec: smooth}.



\subsubsection{A broader conditional distribution class with smoothness disparity}


 In Appendix \ref{sec: approx}, we present a novel approximation result for the function class exhibiting \textit{smoothness disparity} in \Cref{thm.approx_network_One_fct2}.  This new result facilitates the study of theoretical properties of  estimators when the generator \( G_* \in \mathcal{H}_{\dd, \pp}^{\beta_Z, \beta_X} \left(\calZ, \calX, K \right) \).
Note that such a function class defined in \eqref{eq: holder disparity} in Appendix \ref{sec: approx} is much broader compared to the smoothness class in  Section \ref{sec-global} as \( Z \) and \( X \) do not have to be jointly smooth and it  allows for smoothness disparity among them.  The subsequent \Cref{thm: disparity} combines our approximation result with \eqref{eq: entropy Lin} and enables us to specialize \Cref{thm: main} to this class (see \Cref{sec: proof thm: disparity} for the proof).

\begin{theorem}\label{thm: disparity}

Let \( G_* \in \mathcal{H}_{\dd, \pp}^{\beta_Z, \beta_X} \left(\calZ, \calX, K \right) \). Suppose that Assumption \ref{assume: basic} holds and $\sigma_* \in \left[ \sigma_{\min}, \sigma_{\max} \right]$ with $\sigma_{\min} \leq 1$ and $\sigma_{\max} < \infty$. Moreover, we assume $\sigma_* \asymp n^{-\alpha}$, and $\sigma_{\min} = n^{-\gamma}$ for some $ 0 \le \alpha \le \gamma \le (\beta_Z^{-1}\dd + \beta_X^{-1}\pp)^{-1}$. Then, for every $\delta_{\app} \in [0,1]$, we have: Let $\F_s = \F_s\lt(L, r, s, 1, K \rt)$ with $L \asymp \log \delta^{-1}_{\app} $, $r \asymp \delta^{-(\beta_Z^{-1}\dd + \beta_X^{-1}\pp )}_{\app}$, $s \asymp \delta^{-(\beta_Z^{-1}\dd + \beta_X^{-1}\pp )}_{\app} \log \delta^{-1}_{\app}$. Then the sieve MLE $\what{p}_s$ satisfies \eqref{eq: sieve MLE convergence sparse} with the rate outlined in \eqref{eq: rate} with $\xi = \delta^{-(\beta_Z^{-1}\dd + \beta_X^{-1}\pp )}_{\app} \log^2 \delta^{-1}_{\app} $ and $A = \log^2 \delta^{-1}_{\app}$, provided that $\eta_n \le n \eps^{*2}_n/6$.
In particular, choosing 
$
\delta_{\app} := \lt( \sigma_*^2/n \rt)^{1/\lt(2 + \beta_Z^{-1}\dd + \beta_X^{-1}\pp \rt)} \le 1
$
minimizes $\eps_n^* \asymp \sqrt{ \xi \lt\{ A + \log \lt( n/\sigma_{\min} \rt) \rt\}/n} \vee \delta_{\app}/\sigma_* $, and gives
\begin{equation}\label{eq: eps disparity}
    \eps_n^* \asymp n^{-\frac{1 -  \alpha(\beta_Z^{-1}\dd + \beta_X^{-1}\pp)}{2 + \beta_Z^{-1}\dd + \beta_X^{-1}\pp}} \log^{2} (n).    
\end{equation}
\end{theorem}
The proof of \Cref{thm: disparity} is provided in \Cref{sec: proof thm: disparity}.
In the special case when $\alpha = 0$ and $\dd=D$, our convergence rate in \eqref{eq: eps disparity} recovers the minimax optimal rate for conditional density estimation based on kernel smoothing, as established in \cite{li2022minimax}.

\subsubsection{Conditional distribution on manifolds} 
In this part, we extend \Cref{lemma: Caffareli} and provide the existence of the generator when the conditional distribution is supported  on a compact manifold with dimension $\mathsf{d}_\ast \le D$. Due to space constraints, we provide only a sketched proof here; the detailed proof can be found in \Cref{sec: manifold}. Specifically, we first present arguments for the existence of the generator when \( \calY \) is covered by a single chart. We then extend this to the multiple chart case using the technique of partition of unity.

In the simpler case when there exists a single \((\calY, \varphi)\) covering \(\calY\), where \(\varphi: \calB_1(\bzero_{\mathsf{d}_\ast}) \to \calY\) is a homeomorphism, we assume \(\varphi \in \h^{\beta_{\min}+1}\). In this case, we use the change of variable formula to transfer the measure on \(\calB_1(\bzero_{\mathsf{d}_\ast})\) (unit ball in $\real^{\mathsf{d}_\ast}$) from \(\calY\). Following \Cref{lemma: Caffareli}, we can find a transport map \(g \in \h^{\beta_{\min}}\) mapping from \(\calZ \times \calX\) to \(\calB_1(\bzero_{\mathsf{d}_\ast})\). The map \(g \circ \varphi\) then serves as our generator.

In the general case where the compact manifold \(\mathcal{Y}\) needs to be covered by multiple charts, demonstrating the existence of a transport or push-forward map is challenging because \(\mathcal{Y}\) is not uniformly convex. Suppose that \(\{(U_k, \varphi_k)\}_{k=1}^K\) forms a cover of \(\calY\). Due to the compactness of \(\calY\), the number of charts \(K\) is finite. 
Analogous to the single chart scenario, we first construct \(g_k \circ \varphi_k\) to transport the measure on each chart. We then patch these local transport maps together to construct a global transport map; see  \Cref{sec: manifold} for full details. As a result, following \Cref{cor: both case}, the convergence rate for density estimation shall be given by \(\eps^*_n \asymp n^{- (\beta_{\min} - \dd\alpha)/(2\beta_{\min} + \dd)}\).

\section{Numerical Results}\label{sec: Simulation}
In this section, we present numerical experiments to validate and complement our theoretical findings using two synthetic dataset examples. These experiments cover a range of scenarios, including full-dimensional cases as well as benchmark examples involving manifold-based data. Additionally, we provide a real data example to further enrich our experimentation and validation process.
It is worth noting that, although not significant, the computational cost of fitting a conditional generative model is higher compared to fitting an unconditional one, as the input dimension of the deep neural network (DNN) is $\pp + \dd$ rather than just $\dd$.

\medskip
\noindent
\textbf{Learning algorithm to compute sieve MLE.}
\\
For the computational algorithm, we adopt a common conditional variational auto-encoder (VAE) architecture to maximize the following log-likelihood term:$\sum_{j=1}^n \calL_{\mathrm VAE}(g, \sigma, \phi; Y_j, X_j) $, where

\begin{equation*}
    \calL_{\mathrm VAE}(g, \sigma, \phi; y, x) = \log \lt( \frac{p_{g,\sigma}(y,x,z)}{q_{\phi}(Z|y,x)}   \rt).
\end{equation*}

The variational distribution $q_{\phi}(Z|y,x)$ is chosen as the standard normal family $\normal (\mu_\phi(y,x), \Sigma_\phi(y,x))$.

We examine two classes of datasets: (i) full-dimensional response and (ii) response residing on a low-dimensional manifold. The first highlights the generality of our proposed approach, while the second underscores its efficiency in terms of the Wasserstein metric and validates the small noise perturbation strategy outlined in \Cref{cor: perturbation}.

\medskip
\noindent
\textbf{Simulation from full dimension distribution}. We use the following models for data generation.
\begin{itemize}[leftmargin=*,itemsep=-2pt]
    \item \textbf{FD1} : $Y = \mathbb{I}_{\{U<0.5\}}\, \normal \lt(-X, 0.25^2 \rt) + \mathbb{I}_{\{U >0.5\}}\, \normal \lt(X, 0.25^2 \rt)$; $U\sim \mathrm{Unif}(0,1)$, $X \sim \normal(3, 1)$.
    \item \textbf{FD2} : $ Y = X_1^2 + e^{(X_2+X_3/3)} + \sin(X_4 + X_5) + \eps $; $\{X_j\}_{j=1}^5 \stackrel{i.i.d}{\sim} \normal(0,1) $, $\eps \sim \normal(0,1) $.
    \item \textbf{FD3} : $ Y = X_1^2 + e^{(X_2+X_3/3)} + X_4 - X_5 + 0.5\,(1 + X_2^2 + X_5^2) \times \eps $; $\{X_j\}_{j=1}^5 \stackrel{i.i.d}{\sim} \normal(0,1)$, $\eps \sim \normal(0,1) $.
\end{itemize}

These are examples of a mixture model, an additive noise model, and a multiplicative noise model, respectively. The neural architecture for both the encoder and decoder consists of two deep layers, i.e., $L=2$. The hyperparameters are as follows: $r_{\mathrm{enc}} = (\pp+1, 10, 10)$ for $\mu_\phi$ and $\Sigma_\phi$, and $r_{\mathrm{dec}} = (10 + \pp, 10, 1)$ for $g$. The sample size used for simulation is 5000, with a training-to-testing ratio of $4:1$. We employ a batch size of $64$ with a learning rate of $10^{-3}$.

We compare the sieve MLE with CKDE \citep{hall2004cross} and FlexCode proposed by \citet{Flexcode}. To evaluate their performance, we compute the mean squared error (MSE) for both the mean and the standard deviation.
We use Monte Carlo approximation to compute the mean and standard deviation for the sieve MLE, and numerical integration for CKDE and Flexcode. This evaluation strategy resembles that implemented by \citet{huang2022sampling}. \Cref{table:Full_data} summarizes the findings.
\begin{table}[htpb]
\caption{MSE for the estimated conditional mean and the standard deviation.}
\vskip -1em
\begin{center}
\begin{small}
\begin{tabular}{llccc}
\toprule
& & Sieve MLE & CKDE & FlexCode \\
\midrule
FD1 & MEAN & \textbf{0.0379} $\pm$ 0.0170 & 1.0053 $\pm$ 0.1004 & 1.1660 $\pm$ 0.1076 \\
    & SD   & \textbf{0.0280} $\pm$ 0.0045 & 0.9887 $\pm$ 0.0347 & 1.2000 $\pm$ 0.0126 \\
FD2 & MEAN & \textbf{0.1943} $\pm$ 0.0427 & 0.2640 $\pm$ 0.0515 & 0.3954 $\pm$ 0.0571 \\
    & SD   & \textbf{0.2843} $\pm$ 0.0093 & 0.2853 $\pm$ 0.0213 & 5.8278 $\pm$ 0.1607 \\
FD3 & MEAN & \textbf{0.2337} $\pm$ 0.0453 & 0.2967 $\pm$ 0.0537 & 1.3419 $\pm$ 0.1087 \\
    & SD   & 1.6394 $\pm$ 0.0861        & \textbf{0.6334} $\pm$ 0.0460 & 11.4898 $\pm$ 0.1559 \\
\bottomrule
\end{tabular}
\label{table:Full_data}
\end{small}
\end{center}
\vskip -1em
\end{table}

Note that the sieve MLE outperforms all other methods in all scenarios except for the MSE(SD) for the FD3 dataset. However, for the FD3 dataset, we found that as the training sample size increases further, the MSE(SD) of the sieve MLE achieves performance increasingly comparable to CKDE.

\medskip
\noindent
\textbf{Simulation from distributions on manifolds.}
We consider two examples of manifolds with an intrinsic dimension $\dd=1$, while the ambient dimension is $D=2$. 
\begin{itemize}[leftmargin=*,itemsep=0pt]
    \item \textbf{M1} : $Y = G_*(Z,U) + \eps$, with $G_* = (G_*^{(1)}, G_*^{(2)}) $, $G_*^{(1)} = \mathbb{I}_{\{U < 0.5\}} \lt( 1 -  \cos(Z) \rt) + \mathbb{I}_{\{U > 0.5\}} \cos(Z)$, $G_*^{(2)} =  \mathbb{I}_{\{U < 0.5\}} \lt( 0.5 -  \sin(Z) \rt) + \mathbb{I}_{\{U>0.5\}} \sin(Z)$; $Z\sim \mathrm{Unif}(0,\pi)$, $U\sim \mathrm{Unif}(0,1)$.
    \item \textbf{M2} : $Y = G_*(Z,U) + \eps$, with $G_* = \lt(G_*^{(1)}, G_*^{(2)}\rt) $, $G_*^{(1)} = \mathbb{I}_{\{U<0.5\}} \cos(Z) + \mathbb{I}_{\{U > 0.5\}} \, 2\cos(Z) $, $G_*^{(2)} = \mathbb{I}_{\{U<0.5\}} \, 0.5 \sin(Z) + \mathbb{I}_{\{U > 0.5\}}  \sin(Z)$; $Z\sim \mathrm{Unif}(0,2\pi)$, $U\sim \mathrm{Unif}(0,1)$.
\end{itemize}
The manifold $M_1$ consists of two moons. The manifold $M_2$ comprises ellipses, with conditions distinguishing the inner and outer confocal ellipses. The noise sequence follows a two-dimensional centered Gaussian distribution, $\eps \sim \normal(0_2, \sigma_*^2 I_2)$. We investigated this setup across various noise variances $\sigma_*^2$. Our neural architecture employed $r_{\mathrm{enc}} = (\pp+2, 100, 100, 2)$ for $\mu_\phi$ and $\Sigma_\phi$, and $r_{\mathrm{dec}} = (2+ \pp, 100, 100, 2)$ for $g$. We utilized a sample size of 5000 for simulation, with a training-to-testing ratio of $4:1$. A batch size of $100$ was employed, with a learning rate of $10^{-3}$.

\begin{figure}[h]
    \centering
    \begin{subfigure}{0.23\textwidth}
        \centering
        \includegraphics[width=\linewidth, height = 0.9\linewidth]{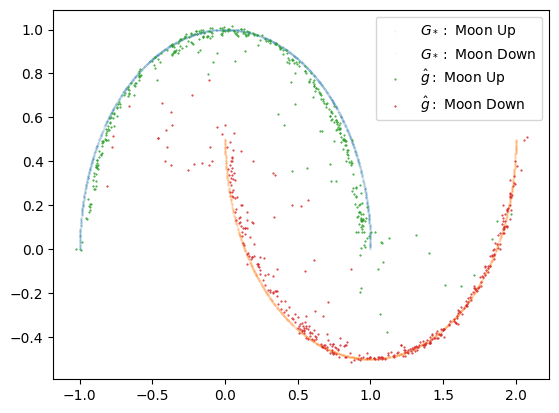}
    \end{subfigure}
    \begin{subfigure}{0.23\textwidth}
        \centering
        \includegraphics[width=\linewidth, height = 0.9\linewidth]{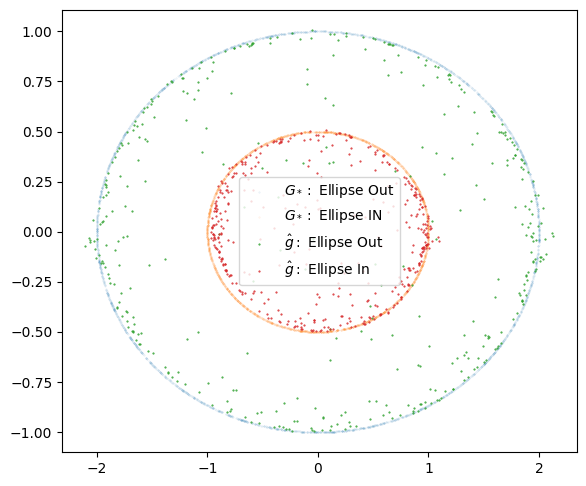}
    \end{subfigure}
    \hfill
    \begin{subfigure}{0.23\textwidth}
        \centering
        \includegraphics[width=\linewidth, height = 0.9\linewidth]{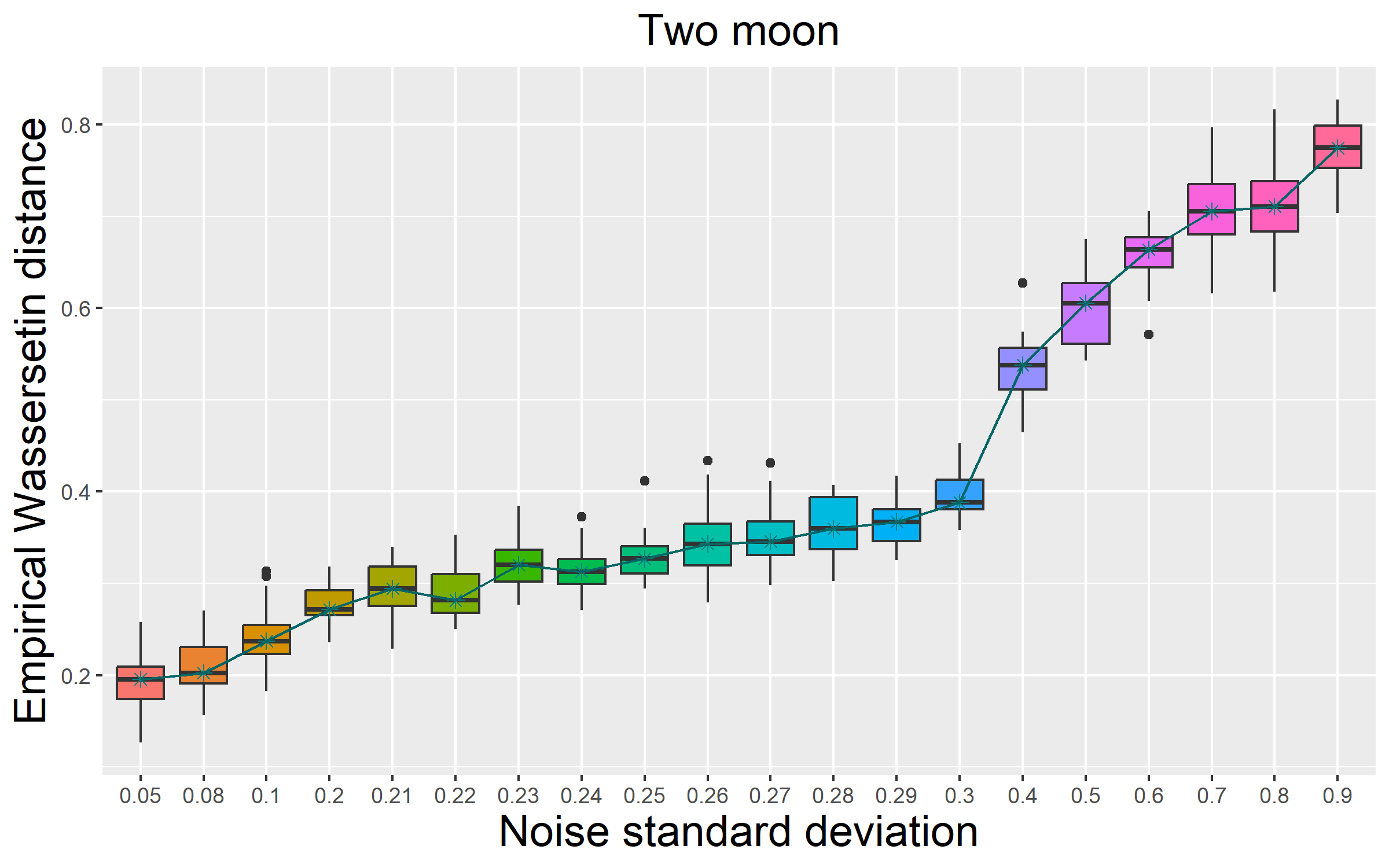}
    \end{subfigure}
    \begin{subfigure}{0.23\textwidth}
        \centering
        \includegraphics[width=\linewidth, height = 0.9\linewidth]{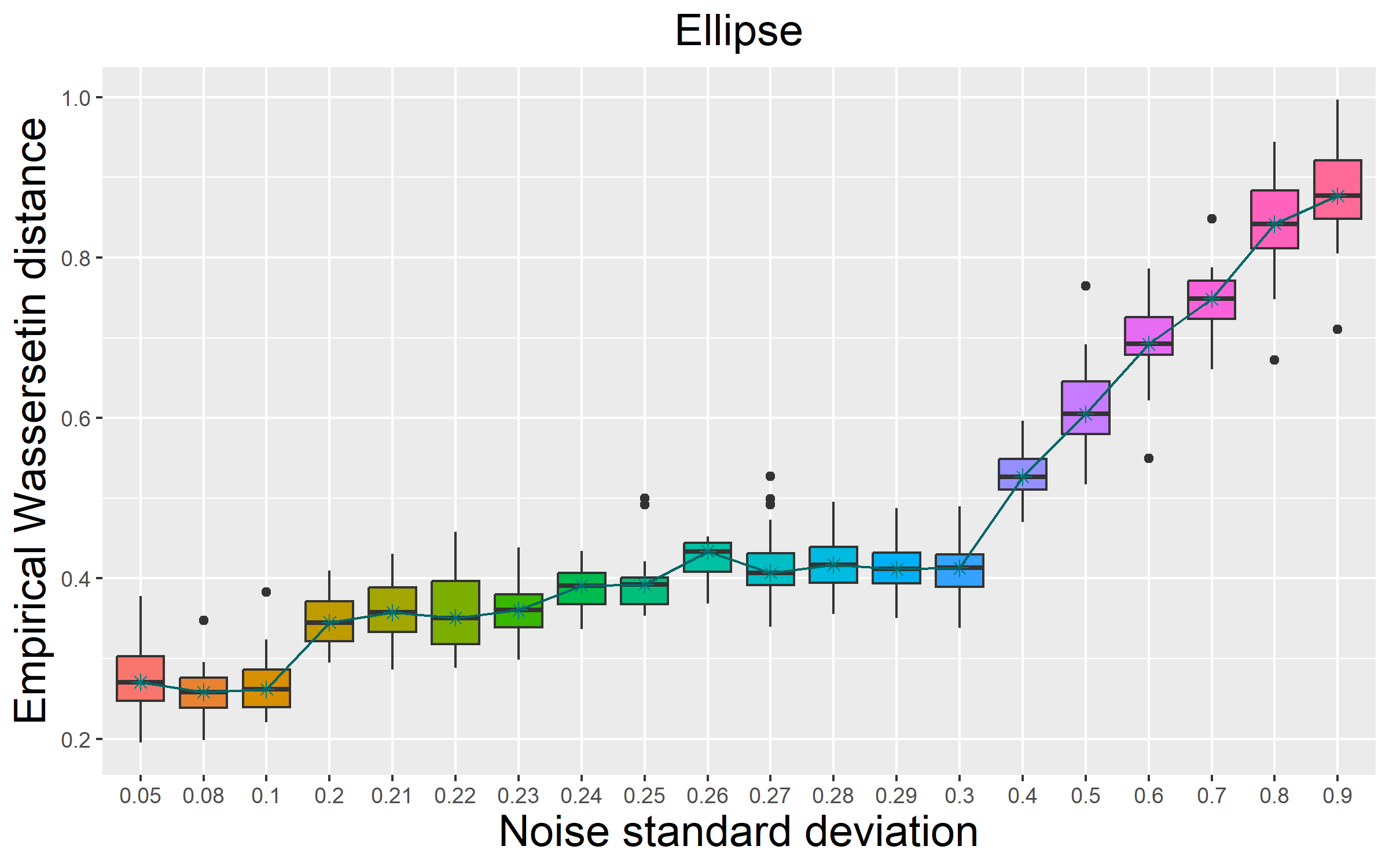}
    \end{subfigure}
    \caption{Generated samples from manifold  $M_1$ and $M_2$ are displayed in the left panel.  The right panel shows box plots for the empirical Wasserstein distance at different noise levels $\sigma_*$.}\label{fig:manifold_images}
\end{figure}


We computed the empirical $W_1$ distance using the algorithm proposed by \citet{cuturi2013sinkhorn} to evaluate the performance. \Cref{fig:manifold_images} (right panel) presents the box plots of $W_1$ between the true and learned distribution for $M_1$ and $M_2$ across 20 repetitions. We observe the following general behaviors:
\begin{itemize}[leftmargin=*,itemsep=0pt]
    \item When $\alpha$ is small and close to zero, the noise variance is large, making estimation challenging due to the singularity of the true data distribution.
    \item When $\alpha$ is large, the noise variance is small, and the perturbed data facilitates efficient estimation.
\end{itemize}
This observed pattern, as emphasized in \Cref{cor: perturbation}, closely aligns with the results achieved in \eqref{eq: manifold rate}. An additional numerical experiment on real data has been performed and can be found in \Cref{sec: Real data}.

\section{Discussion}
 We investigated  statistical properties of a likelihood-based conditional deep generative model for distribution regression in a scenario where the response variable is situated in a high-dimensional ambient space but is centered around a potentially lower-dimensional intrinsic structure. Our analysis established favorable rates in both the Hellinger and Wasserstein metrics which are dependent on only the intrinsic dimension of the data. Our theoretical findings show that the conditional deep generative models can circumvent the curse of dimensionality for high-dimensional distribution regression. To the best of our knowledge, our work is the first of its kind.

Given the novelty of emerging statistical methodologies with intricate structural considerations in the study of deep generative models, there exist numerous paths for future exploration. Among these potential directions, we are particularly interested in investigating controllable generation via penalized optimization methods, studying statistical properties of deep generative models trained via matching flows, as well as delving into the hypothesis testing problem within the framework of deep generative models, among others. Another interesting direction is to explore residual neural network structure for modeling time series of distributions with interesting temporal dependence structures.


\bibliography{references}
\bibliographystyle{apalike}

\newpage
\appendix
\onecolumn

\begin{center}
    {\bf \Large 
Supplementary Materials for ``A Likelihood Based Approach to Distribution Regression Using Conditional Deep Generative Models''
 }
\end{center}


\section{Additional numerical results}
\subsection{Numerical result for real data}\label{sec: Real data}
We used the MNIST dataset for two purposes: to demonstrate the generalizability of our approach to a benchmark image dataset where the intrinsic dimension $\dd$ is much lesser than the ambient dimension $D=784$ and to underscore the effectiveness of sparse networks as outlined in \Cref{lemma: DNN approximation}.1 and \Cref{cor: both case}.1.

For the fully connected architecture, we set $r_{\mathrm{enc}} = (10+784, 512,  2)$ for $\mu_\phi$ and $\Sigma_\phi$, and $r_{\mathrm{dec}} = (10 + 2, 512, 784)$ for $g$. For the sparse architecture, we use $r_{\mathrm{enc}} = (10+784, 608, 432, 256,  2)$ for $\mu_\phi$ and $\Sigma_\phi$, and $r_{\mathrm{dec}} = (10 + 2, 256, 432, 608, 784)$ for $g$. The input dimension of $10$ for both the encoder and decoder corresponds to the one-hot encoding of the labels. We employ a batch size of $64$ with a learning rate of $10^{-3}$.

\begin{figure}[H]
    \centering
    \begin{subfigure}{0.3\textwidth}
        \centering
        \includegraphics[width=\linewidth]{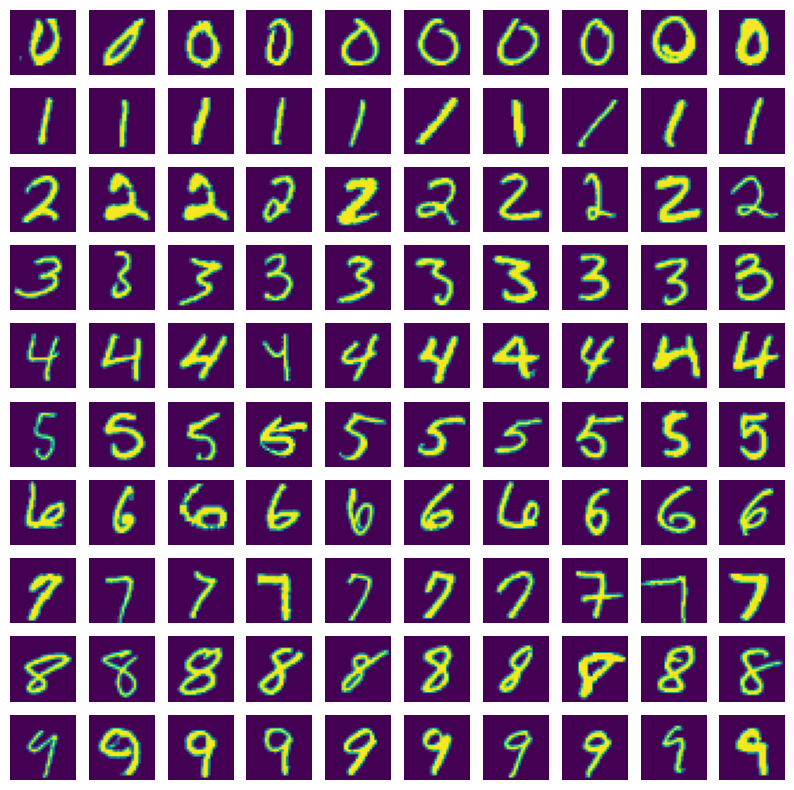}
    \end{subfigure}
    \begin{subfigure}{0.3\textwidth}
        \centering
        \includegraphics[width=\linewidth]{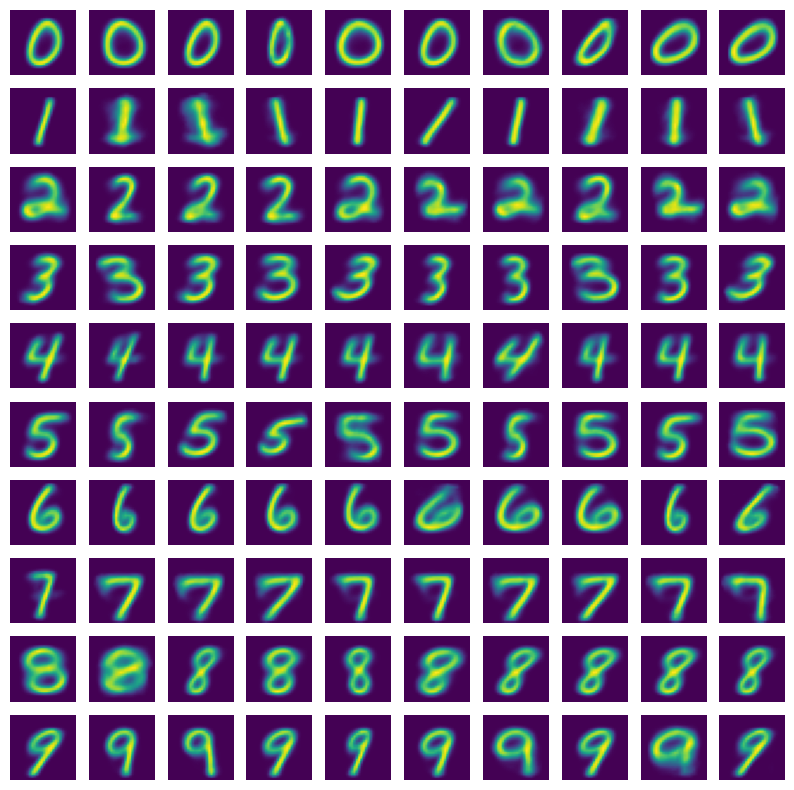}
    \end{subfigure}
    \begin{subfigure}{0.3\textwidth}
        \centering
        \includegraphics[width=\linewidth]{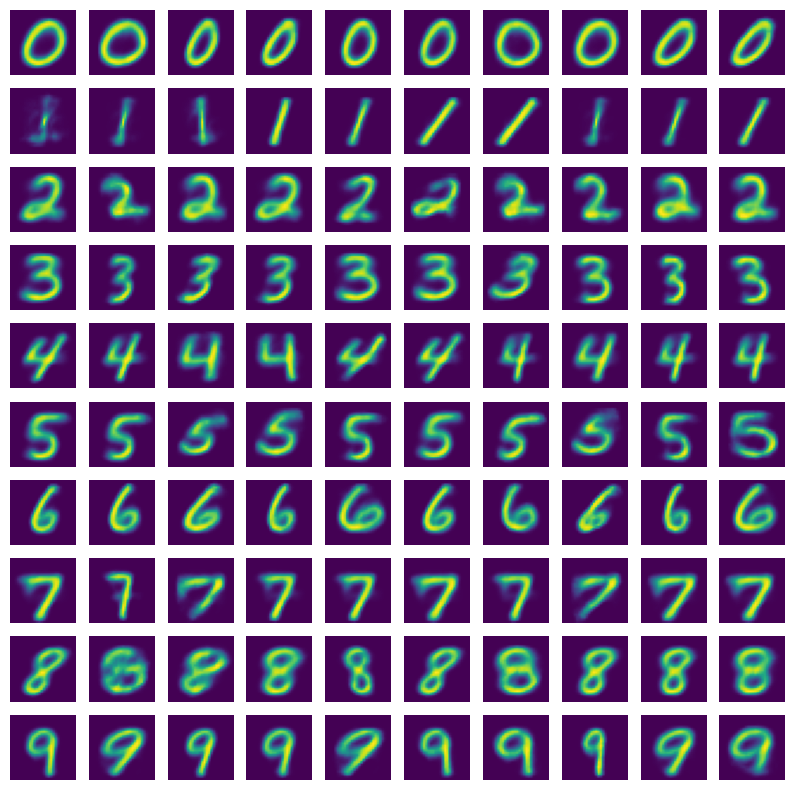}
    \end{subfigure}
    \caption{MNIST images: real images (left panel), generated images with sparse architecture (central panel), and generated images with fully connected architecture (right panel)}
    \label{fig:mnist_images}
\end{figure}

\Cref{fig:mnist_images} presents a visual comparison between real and generated images, organized according to their respective labels. The real images were randomly sampled from the training set along with their corresponding labels, while the generated images were produced using these labels (conditions) and random seeds.

This MNIST example highlights a case where the intrinsic dimension is significantly smaller than the ambient data dimension. This example serves to validate the proposed methodology in high-dimensional settings.

To quantify sample quality, we computed the Wasserstein-1 distance ($W_1$) between generated and test images. For each digit, we averaged $W_1$ distances over 50 samples, reporting results as mean $\pm$ standard deviation. For reference, the baseline $W_1$ distance between two test images is $2.0219 \pm 0.7450$. \Cref{table:Full_data2} summarizes these distances across different levels of Gaussian noise added during training.
\begin{table}[H]
    \centering
    \caption{Mean W1 distance ($\pm$ SD) between generated and test MNIST images under varying training-data noise.}
    \label{table:Full_data2}
    \begin{tabular}{l@{\hspace{1cm}}c@{\hspace{1cm}}c}
    \toprule
    Noise & Sparsely connected & Fully connected \\
    \midrule
    0     & $1.9555 \pm 0.7182$ & $1.8859 \pm 0.7355$ \\
    0.005 & $1.9478 \pm 0.7329$ & $1.8663 \pm 0.6251$ \\
    0.01  & $1.9503 \pm 0.7291$ & $1.9598 \pm 0.6867$ \\
    0.02  & $2.0699 \pm 0.6937$ & $2.0616 \pm 0.7410$ \\
    0.04  & $2.2199 \pm 0.6735$ & $2.2117 \pm 0.6627$ \\
    0.06  & $2.3487 \pm 0.6576$ & $2.3172 \pm 0.6267$ \\
    0.08  & $2.4623 \pm 0.6245$ & $2.4076 \pm 0.6308$ \\
    0.1   & $2.5734 \pm 0.6492$ & $2.5002 \pm 0.6337$ \\
    0.3   & $3.4931 \pm 0.7012$ & $3.4943 \pm 0.7164$ \\
    0.5   & $4.0880 \pm 0.7518$ & $4.0995 \pm 0.7633$ \\
    \bottomrule
    \end{tabular}
\end{table}
As shown in \Cref{table:Full_data2}, at zero noise both architectures achieve W1 distances slightly below the baseline, indicating high-fidelity sample generation. As noise increases, W1 distances grow steadily, reflecting degradation in sample quality. Both network types follow similar trends, underlining robustness to architectural choice; minor deviations suggest subtle differences in sensitivity to noise. These empirical observations accord with our theoretical predictions on the large-sample properties of manifold-supported data.

\subsection{Additional numerical results for distributions on manifold}
We extended our analysis to examine how the empirical $W_1$ distance varies with sample size, while keeping the noise level fixed at $\sigma_* = 0.01$. Below is a summary table showing the median empirical Wasserstein distances for different sample sizes. The experimental setup remains consistent with the manifold case described in the \Cref{sec: Simulation}.

\begin{table}[H]
    \centering
    \caption{Empirical Wasserstein distance $W_1$ (median) for different sample sizes}
    \begin{tabular}{ccc}
    \toprule
    Sample Size & Two Moon ($\sigma_* = 0.01$) & Ellipse ($\sigma_* = 0.01$) \\
    \toprule
    4000 & 0.251 & 0.295 \\
    6000 & 0.232 & 0.285 \\
    7000 & 0.216 & 0.271 \\
    8000 & 0.214 & 0.253 \\
    9000 & 0.212 & 0.259 \\
    10000 & 0.196 & 0.251 \\
    \bottomrule
    \end{tabular}
\end{table}

While extracting exact rates through simulation can be challenging, the results in the table validate the large-sample properties for manifolds. These empirical findings align well with the theoretical expectations, further confirming the consistency and convergence trends of our framework.

\section{Notation}\label{sec: notations}
We denote \( a \vee b \) and \( a \wedge b \) as the maximum and minimum of two real numbers \( a \) and \( b \), respectively. The notation \( \lceil a \rceil \) represents the smallest integer greater than or equal to \( a \). The inequality \( a \lesssim b \) indicates that \( a \) is less than or equal to \( b \) up to a multiplicative constant. When we write \( a \lesssim_{\log} b \), it means that \( a \) is less than or equal to \( b \) up to a logarithmic factor, specifically \(\log(n)\). We denote \( a \asymp b \) when both \( a \lesssim b \) and \( b \lesssim a \) hold. For vector norms, \( |\cdot|_{\pp} \) represents the \(\ell^{\pp}\) norm, while \(\|\cdot\|_{\pp}\) denotes the \(L^{\pp}\)-norm of a function for \( 1 \le \pp \le \infty \). Lastly, \(\calB_{\epsilon}(u)\) signifies the Euclidean open ball with radius \(\epsilon\) centered at \(u\).

We use the multi-index notation through the main paper and the appendix. Denote $\mathbb N$ as the set of natural numbers and  $\nat$ as $\mathbb{N}\cup \{0\}$. For a vector $\bx \in \real^r$, we denote the components as $\bx = (x^\1, \ldots, x^\upr) $. Given a function $f: D \subset \real^r \to \real$, the operator is defined as $\dau^{\balpha} := \dau^{\alpha^\1} \ldots \dau^{\alpha^\upr}$ with $\balpha \in \nat^r$, where $\dau^{\alpha^{(j)}} f := \dau^{\alpha^{(j)}} f(\bx)/\dau x^{(j)}$. For $\balpha \in \nat^r$, the expression $|\balpha| = \sum_{j=1}^r |\alpha^{(j)}|$. Given a function $f(\cdot, \cdot) : D \times D_\prime \subset \real^{r}\times \real^{r_\prime} \to \real$, we denote the operator $\dau^{\balpha+\balpha_\prime} := \dau^{\alpha^\1} \ldots \dau^{\alpha^{(r)}}\,\dau^{\alpha_\prime^\1} \ldots \dau^{\alpha_\prime^{(r_\prime
)}}$, with $\balpha \in \nat^r$ and $\balpha_\prime \in \nat^{r_\prime}$, where $\dau^{\alpha^{(j)}} f(\bx,\by) = \dau^{\alpha^{(j)}} f(\bx,\by)/\dau^{\alpha^{(j)}} x^{(j)} $ and $\dau^{\alpha_\prime^{(j)}} f(\bx,\by) = \dau^{\alpha_\prime^{(j)}} f(\bx, \by)/\dau y^{(j)}$, with $\bx \in D$ and $\by \in D_\prime$. This notation allows us to represent the derivative with variable $\bx$ and $\by$ separately through the vector $\balpha$ and $\balpha_\prime$ which is required to tackle the smoothness disparity along $x$ and $y$ variable. The $\beta-$H\"older class functions are defined as
\begin{equation}\label{eq: holder class}
    \begin{aligned}
        \mathcal{H}_{{r}}^{\beta} &\left( D, M \right) = \Big\{ 
	       f :D \subset \real^r \rightarrow \mathbb{R} : \\
	       &\sum_{ \substack{{\balpha} : |{\balpha}| < \beta } }\|\partial^{{\balpha} } f\|_\infty + \sum_{ \substack{{\balpha} : |{\balpha} |= \lfloor \beta \rfloor }} \, \sup_{\substack{\bu_1, \bu_2 \in D \\  \bu_1 \neq \bu_2}}
	       \frac{|\partial^{{\balpha}  } f(\bu_1) - \partial^{{\balpha} } f( \bu_2)|}{|\bu_1 - \bu_2|_\infty^{\beta-\lfloor \beta \rfloor}} \leq M
	       \Big\},
    \end{aligned}
\end{equation}
We extend this definition to include the H\"older class of functions with differences in smoothness (smoothness disparity) along two variables. This class is defined as
\begin{equation}\label{eq: holder disparity}
    \begin{aligned}
        \mathcal{H}_{{r}, {r_\prime}}^{\beta, \beta_\prime} &\left( D, D_\prime, M \right) = \Big\{ 
    	f(\cdot, \cdot) : D \times D_\prime \subset \mathbb{R}^{{r}} \times \mathbb{R}^{{r_\prime}} \rightarrow \mathbb{R} : \\
    	&\sum_{ \substack{{\balpha} : |{\balpha}| < \beta \\ {\balpha_\prime} : |{\balpha_\prime}| < \beta_\prime} }\|\partial^{{\balpha} + {\balpha_\prime}} f\|_\infty + \sum_{ \substack{{\balpha} : |{\balpha} |= \lfloor \beta \rfloor \\ {\balpha_\prime} : |{\balpha_\prime} |= \lfloor \beta_\prime \rfloor } } \, \sup_{\substack{\bu_1, \bu_2 \in D_X \\ \bv _1, \bv _2 \in D_Y \\  \bu_1 \neq \bu_2 \\ \bv _1 \neq \bv _2 }}
    	\frac{|\partial^{{\balpha} + {\balpha_\prime}} f(\bv _1, \bu_1) - \partial^{{\balpha} + {\balpha_\prime}} f(\bv _2, \bu_2)|}{|\bu_1 - \bu_2|_\infty^{\beta-\lfloor \beta \rfloor} \vee |\bv _1 - \bv _2|_\infty^{\beta_\prime-\lfloor \beta_\prime \rfloor}} \leq M
    	\Big\}.
    \end{aligned}
\end{equation}
We denote $\mathcal{H}_{{r}}^{\beta} (D) = \cup_{M>0} \mathcal{H}_{{r}}^{\beta} \left( D, M \right) $ and $\mathcal{H}_{{r}, {r_\prime}}^{\beta, \beta_\prime} \left( D, D_\prime \right) = \cup_{M>0} \mathcal{H}_{{r}, {r_\prime}}^{\beta, \beta_\prime}\left( D, D_\prime, M \right)$.

\section{More on Smooth conditional density}\label{sec: smooth}
\begin{theorem}[\cite{villani2009optimal} Theorem 12.50]\label{lemma: transport}
Suppose that
    \begin{enumerate}[label=(\roman*)]
        \item $\mathcal{A}_1 $ and $\mathcal{A}_2 $ are uniformly convex, bounded, open subsets of $\real^\dd$ with $\mathcal{C}^{\lfloor \beta \rfloor + 2}$ (continuously differentiable up to order ${\lfloor \beta \rfloor + 2}$) boundaries,
        \item $h_1 \in \h^\beta(\mathcal{A}_1)$ and $h_2 \in \h^\beta(\mathcal{A}_2)$ for some $\beta>0$, are probability densities bounded above and below.
    \end{enumerate}
        Then, there exists a unique map (up to an additive constant) $g: \mathcal{A}_1 \to \mathcal{A}_2$ with $g \in \h^{\beta+1}(\mathcal{A}_1)$, such that if $U\sim h_1$ then $g(U) \sim h_2$.
\end{theorem}

\begin{proof}[Proof of \Cref{lemma: Caffareli}]
    Given that $Z$ and $X$ is independent, the product measure on $\calZ \times \calX$ is $p_Z \mu_X^*$. Following the smoothness from $p_Z$ and $\mu_X^*$, the map $p_Z(\cdot) \mu_X^*(\cdot) \in \h^{\min\{\beta_Z, \beta_X\}}(\calZ \times \calX)$. This implies that $p_Z(\cdot) \mu_X^*(\cdot) \in \h^{\min\{\beta_Z, \beta_X, \beta_Q\}}(\calZ \times \calX)$. Again $q_* \in \h^{\beta_Q}(\calY)$ implies $q_* \in \h^{\min\{\beta_Z, \beta_X, \beta_Q\}}(\calY)$. The result now follows directly from \Cref{lemma: transport}.
\end{proof}

Many of the problems in the conditional setting have an analog in the joint setup. Our proposed approach has a direct statistical extension to this setup. The sufficiency of such extension follows from the observation in the subsequent \Cref{lemma: existence} which is based on Lemma 2.1 and Lemma 2.2 of \citet{huang2022sampling} (see also Theorem 5.10 of \citet{kallenberg1997foundations}).

\begin{lemma}[Noise Outsourcing Lemma]\label{lemma: existence}
    Let $(Y,X) \in \calY \times \calX$ with joint distribution $P_{Y,X}$. Suppose $Y$ is standard Borel space, then there exists $Z \sim \normal(0, I_m)$ for any given $m\ge 1$, independent of $X$, and a Borel measurable function $G: \R^m \times \calX \to \calY $ such that
    \begin{equation}\label{eq: existence of G}
    \lt(X, G(Z,X)\rt) \sim (Y,X).    
    \end{equation}
    Moreover, the condition \eqref{eq: existence of G} is equivalent of 
    $$
    G(Z,x) \sim P_{Y| X=x}.
    $$
\end{lemma}

\section{More on Conditional distribution on manifolds}\label{sec: manifold}

Suppose \((\calY, \varphi)\) is the single chart covering \(\calY\), where \(\varphi: \calB_1(\bzero_{\mathsf{d}_\ast}) \to \calY\) is a homeomorphism. We assume that \(\varphi \in \h^{\beta_{\min}+1}\), and that $\inf_{\bu \in \calB_1(\bzero_{\mathsf{d}_\ast})} |J_{\varphi}(\bu)|$ is bounded below by a positive constant, where
$$
|J_{\varphi}(\bu)| = \sqrt{\mathrm{det} \lt( \frac{\dau \varphi}{\dau \bu^\top} \frac{\dau \varphi}{\dau \bu} \rt)}
$$
is the Jacobian determinant of $\varphi$.

Note that when $\mathsf{d}_\ast < D$, the distribution $Q_*$ cannot possess a Lebesgue density because of the singularity of $\calY$. We, therefore consider a density with respect to the $\mathsf{d}_\ast-$dimensional Hausdorff measure in $\real^D$, denoted by $\mathsf{H}_{\mathsf{d}_\ast}$. Suppose that $Q$ allows the Radon-Nikodym derivative $q$ with respect to $\mathsf{H}_{\mathsf{d}_\ast}$. We further assume that $q$ is bounded from above and below and that $q\circ\varphi \in \h^{\beta_{\min}}$. Then by change of variable formula, the Lebesgue density of $\wtilde{Q}$, the push-forward measure on $\calB_1(\bzero_{\mathsf{d}_\ast})$ through the map $\varphi^{-1}$, is given as 
$$
\wtilde{q}(\bu) = q(\varphi(\bu))|J_{\varphi}(\bu)|.
$$
Following the assumptions on the Jacobian determinant and $\varphi \in \h^{\beta_{\min}+1}$, it follows that $|J_{\varphi}(\bu)|$ is bounded from above and below, and the map $\bu \mapsto |J_{\varphi}(\bu)|$ belongs to $\h^{\beta_{\min}}$. Therefore, $\wtilde{q}$ is bounded above and below, belongs to $\h^{\beta_{\min}}(\calB_1(\bzero_{\mathsf{d}_\ast}))$. By \Cref{lemma: Caffareli}, assuming $\beta_{\min} \le \beta_Z \wedge \beta_X $, there exists $g \in \h^{\beta_{{\min}}+1}$ such that $\wtilde{Q} = Q_g$. Thus, we have $Q = Q_{\varphi \circ g}$, where $\varphi \circ g: \calZ \times \calX \to \calY$. Following \Cref{lemma: DNN approximation}, it is possible to find the appropriate neural network approximating them.

Suppose $\calY$ is covered by the charts $\lt\{\lt(U_k, \varphi_k\rt) \rt\}_{k=1}^K$, with $1< K < \infty$, where $\varphi_k : \calB_1(0_{\mathsf{d}_\ast}) \to U_k $ is a homeomorphism. As before, we assume \(\varphi_k \in \h^{\beta_{\min}+1}\), $|J_{\varphi_k}(\bu)|$ is bounded below by a positive constant, $Q$ possesses density $q$ with respect to $\mathsf{H}_{\mathsf{d}_\ast}$ that is bounded above and below, and that $q\circ \phi_k \in \h^{\beta_{\min}}$. Let $Q_k (\cdot) = Q(\cdot)/Q(U_k)$ be the normalized measure of $Q$ over $U_k$.

We denote \( q_k \) as the corresponding density with respect to \(\mathsf{H}_{\mathsf{d}_\ast}\). For \(\bu \in U_k \cap U_\ell\), \( q_k(\bu) Q(U_k) = q_\ell(\bu) Q(U_\ell) = q(\bu) \) holds due to the measure \( Q(\cdot) \) being compatible with the charts. This is ensured because the densities \( Q(U_k) q_k(\cdot) \) and \( Q(U_\ell) q_\ell(\cdot) \) are consistent and align with the measure \( Q \) over the overlapping regions of the charts. This compatibility is essential for constructing a coherent global measure from local chart densities.



A compact manifold \(\calY\) can be covered by a finite partition of unity \(\{\tau_k, k=1, \ldots, K\}\), each sufficiently smooth \citep{lee2012smooth}. By definition, each function in this partition satisfies \(\tau_k(\bu) = 0\) for \(\bu \notin U_k\) and \(\sum_{k=1}^K \tau_k(\bu) = 1\) for all \(\bu \in \calY\). Given that \(q(\bu) = Q(U_k) q_k(\bu)\) for each \(k\) and \(\bu \in U_k\), we can express \(q(\bu)\) as:
\[ q(\bu) = \sum_{k=1}^K Q(U_k) \tau_k(\bu) q_k(\bu). \]
To normalize, let \(c_k = \int \tau_k (\bu) dQ_k(\bu)\) and define \(q^\prime_k(\bu) = \tau_k(\bu) q_k(\bu) / c_k\). Thus, we can rewrite \(q(\bu)\) as:
\[ q(\bu) = \sum_{k=1}^K \pi_k q^\prime_k(\bu), \]
where \(\pi_k = c_k Q(U_k)\). This formulation reveals that \(q\) is a mixture of the component densities \(q^\prime_k(\bu)\), weighted by \(\pi_k\). This mixture approach ensures compatibility across different charts, providing a unified density representation over the entire manifold \(\calY\).

Since \( q^\prime_k \) is sufficiently smooth, we can construct a mapping \( g_k : \widetilde{\mathcal{V}} \to \mathcal{Y} \) such that \( Q^\prime_k \) is the distribution of \( g_k(\widetilde{V}) \), supported on \( U_k \), where \( \widetilde{\mathcal{V}} \) is a uniformly convex set in \( \mathbb{R}^{d_*} \), and \( \widetilde{V} \) follows a uniform distribution on \( \widetilde{\mathcal{V}} \). Next, construct a disjoint partition of the interval \( (0,1) \) into \( K \) intervals \( I_1, \ldots, I_K \) with lengths \( \pi_1, \ldots, \pi_K \), where \( I_k = [\sum_{i=1}^{k-1}\pi_i, \sum_{i=1}^k \pi_i] \). Define \( h_k \) as the indicator function on the interval \( I_k \), i.e., \( h_k(u) = 1 \) if \( u \in I_k \) and \( 0 \) otherwise. For a random variable \( \mathsf{U} \) following Uniform\((0,1)\), it follows that \( P_{\mathsf{U}}(h_k(\mathsf{U}) = 1) = \pi_k \), and \( P_{\mathsf{U}}(h_k(\mathsf{U}) = 0) = 1 - \pi_k \). Now, define \( \mathbf{v} = (\mathsf{u}, \wtilde{v}) \), where \( \mathsf{u} \sim \text{Uniform}(0,1) \) and \( v \sim \text{Uniform}(\widetilde{\mathcal{V}}) \). Using this, construct \( g(\mathbf{v}) = \sum_{k=1}^K h_k(\mathsf{u}) g_k(v) \). It is straightforward to observe that \( Q = Q_g \), as the partitioning through \( h_k \) ensures that the measure is correctly matched to each \( g_k \), and \( g_k \) ensures that the restricted distributions \( Q^\prime_k \) are appropriately supported on \( U_k \).

From an approximation perspective, the indicator functions \( h_k \) and the localized generators can be effectively approximated using ReLU neural networks. This also holds for their products and further linear combinations. For details on such constructions, one may refer to \citet{schmidt2019manifold} for sparse neural networks and \citet{kohler2023manifold} for dense neural networks.

It is important to note that we do not guarantee the regularity of the \( g_k \) maps, as they are not necessarily lower bounded. However, the partition of unity maps \( \tau_k \) vanish only at the boundary of \( U_k \). This property may allows for the construction of sufficiently smooth maps.  For the multiple-chart case, we rely on more stringent results, such as Brenier's Theorem (see, for example, \citet{villani2009optimal}) or the Noise Outsourcing Lemma (\Cref{lemma: existence}), to ensure the existence of the transport maps.

\section{Proof of Lemma~\ref{lemma: Entropy bound}}\label{sec: proof lemma: Entropy bound}
\begin{proof}
For $g_1(\cdot|x), g_2(\cdot|x) \in \F$ with $\| |g_1 - g_2|_\infty\|_\infty \le \eta_1$. Then 
\begin{align}\nonumber
    &p_{g_1, \sigma} (y|x) - p_{g_2, \sigma} (y|x) 
    \\\nonumber
    =& \int \phi_\sigma (y - g_1(x,z)) \lt( 1 - \frac{\phi_\sigma (y - g_2(x,z))}{\phi_\sigma (y - g_1(x,z))} \rt) d P_Z(z)
    \\\nonumber
    =& \int \phi_\sigma (y - g_1(x,z)) \lt( 1 - \exp{\lt\{ - \frac{|y - g_2(x,z)|^2_2 - |y - g_1(x,z)|^2_2}{2\sigma^2} \rt\}} \rt) d P_Z(z)
    \\\label{eq: t1}
    \le& \int \phi_\sigma (y - g_1(x,z)) \lt( \frac{|y - g_2(x,z)|^2_2 - |y - g_1(x,z)|^2_2}{2\sigma^2} \rt) d P_Z(z)
    \\\nonumber
    =& \int \phi_\sigma (y - g_1(x,z)) \lt( \frac{|g_2(x,z) - g_1(x,z)|^2_2 - 2(y - g_1(x,z))^T(g_2(x,z) - g_1(x,z))}{2\sigma^2} \rt) d P_Z(z)
    \\\nonumber
    \le& \int \phi_\sigma (y - g_1(x,z)) \lt( \frac{|g_2(x,z) - g_1(x,z)|^2_2}{2\sigma^2} + \frac{2|y - g_1(x,z)|_1|g_2(x,z) - g_1(x,z)|_\infty}{2\sigma^2} \rt) d P_Z(z)
    \\\label{eq: t2}
    \le& \int \phi_\sigma (y - g_1(x,z)) \frac{2KD \eta_1}{2\sigma^2} d P_Z(z) + \frac{2\eta_1}{2\sigma^2}\int |y - g_1(x,z)|_1 \phi_\sigma (y - g_1(x,z))    d P_Z(z)
    \\\label{eq: t3}
    \le&  \frac{2KD \eta_1}{2\sigma^2} \frac{1}{\lt(\sqrt{2\pi \sigma^2}\rt)^D} + \frac{\eta_1}{\sigma^2} \int  \sqrt{\frac{D}{2\pi e}} \frac{1}{(\sqrt{2\pi \sigma^2})^{D-1}} d P_Z(z)
    \\\label{eq: A1}
    \le& c_1(K, D) \sigma_{\min}^{-(D+2)} \eta_1.
\end{align}
For the last line, we use the fact that $\sigma_{\min} \le 1$. The inequality at \eqref{eq: t1} follows from $e^{-x} \ge (1 - x)$. The ones at \eqref{eq: t2} follows using 
\begin{align*}
    |g_2(x,z) - g_1(x,z)|_2^2  \le 2K|g_2(x,z) - g_1(x,z)|_1  &\le 2K D |g_2(x,z) - g_1(x,z)|_\infty
    \\
    &\le 2KD \||g_1 - g_2|_\infty\|_\infty  \le 2KD\eta_1
\end{align*}
and $|g_2(x,z)- g_1(x,z)|_\infty \le \eta_1$. The change at \eqref{eq: t3} follows from $\phi_\sigma (y - g_1(x,z)) \le \lt(\sqrt{2\pi \sigma^2}\rt)^{-D}$ and the bound
$$
|v|_1\phi_\sigma (v) \le \sqrt{\frac{D}{2\pi e}} \frac{1}{(\sqrt{2\pi \sigma^2})^{D-1}}.
$$

\
\\
Now for $\sigma_1, \sigma_2 \in \lt[ \sigma_{\min}, \sigma_{\max} \rt]$ with $|\sigma_1 - \sigma_2 | \le \eta_2$. It holds that $\lt| \sigma_1^{-2} - \sigma_2^{-2} \rt| \le \sigma_1^{-2} \sigma_2^{-2} \lt( \sigma_1 + \sigma_2\rt) \eta_2$ and $\lt|\log \lt( \frac{\sigma_2}{\sigma_1} \rt) \rt|\le \frac{\eta_2}{\min\{\sigma_1, \sigma_2\}}$. We have
\begin{align}\nonumber
     &p_{g, \sigma_1} (y|x) - p_{g_2, \sigma_2} (y|x)
    \\\nonumber
    =&\int \phi_{\sigma_1}(y - g(x,z) \lt(1 - \lt(\frac{\sigma_1}{\sigma_2}\rt)^D  \exp\lt\{ \frac{|y - g(x,z)|_2^2}{2} \lt( \frac{1}{\sigma_1^2} - \frac{1}{\sigma_2^2} \rt) \rt\} \rt) dP_Z(z)
    \\\label{eq: t4}
    \le & \int \phi_{\sigma_1}(y - g(x,z) \lt[  \frac{|y - g(x,z)|_2^2}{2} \lt( \frac{1}{\sigma_2^2} - \frac{1}{\sigma_1^2} \rt) - D \log \lt(\frac{\sigma_1}{\sigma_2}\rt) \rt] dP_Z(z)
    \\\nonumber
    \le & \int \phi_{\sigma_1}(y - g(x,z) \lt[  \frac{|y - g(x,z)|_2^2}{2} \lt( \frac{\sigma_1 + \sigma_2}{\sigma_1^2\sigma_2^2} \rt) \eta_2 + \frac{D\eta_2}{\min\{\sigma_1, \sigma_2 \}} \rt] dP_Z(z)
    \\\label{eq: t5}
    \le &  \frac{1}{(\sqrt{2\pi \sigma_1^2})^D} \frac{\sigma_1 + \sigma_2}{e\sigma_2^2} \eta_2 + \frac{1}{\lt( \sqrt{2 \pi \sigma_1^2} \rt)^D} \frac{D\eta_2}{\min\{\sigma_1, \sigma_2 \}}
    \\\label{eq: A2}
    \le & c_2(D) \sigma_{\min}^{-(D+1)} \eta_2.
\end{align}
The \eqref{eq: t4} follows from $1 - e^{-\alpha} \le \alpha$. The change at \eqref{eq: t5} follows from $\phi_{\sigma_1}(y - g(x,z)) \le \lt( \sqrt{2\pi \sigma_1^2} \rt)^{-D}$ and
$$
|v|_2^2 \phi_\sigma(v) \le \frac{\sigma^2}{(\sqrt{2\pi \sigma^2})^D}\frac{2}{e}.
$$
Let $\eps > 0$. Let $\{g_1, \ldots, g_{N_1}\}$ be $\eta_1-$covering of $\F$ and $\{\sigma_1, \ldots, \sigma_{N_2}\}$ be $\eta_2-$covering of $\lt[ \sigma_{\min}, \sigma_{\max} \rt]$ with respect to $\| |\cdot|_\infty \|_\infty$ and $|\cdot|_\infty$. By \eqref{eq: A1} and \eqref{eq: A2}, $\eta_1 = c_1^{-1} \sigma_{\min}^{D+2} \eps/4$ and $\eta_2 =  c_2^{-2} \sigma_{\min}^{D+1} \eps/4$ implies
$$
\lt\{ P_{g_i, \sigma_j}(\cdot|\cdot): i = 1, \ldots, N_1, j= 1, \ldots , N_2 \rt\}
$$
forms an $\eps/2-$covering for $\calP$ with respect to $\| \cdot \|_\infty$. Denote the envelope function of $\F$
\begin{align*}
    H(y,x) = \sup_{p \in \calP} p(y|x) &\le \frac{1}{\lt( 2\pi \sigma^2_{\min} \rt)^{-D/2}} \exp\lt\{ -\frac{|y|_2^2 - 4K^2D}{4\sigma^2_{\max}} \rt\} 
    \\
    &= e^{K^2 D/2\sigma_{\max}^2} 2^{D/2} \lt(\frac{\sigma_{\max}}{\sigma_{\min}}\rt)^D \phi_{\sqrt{2}\sigma_{\max}}(y).
\end{align*}
Following from $\int_{|y|_\infty > t} \phi_{\sigma} (y) dy \le 2D e^{-t^2/2\sigma^2}$, we have 
\begin{equation}\nonumber
    \int\int_{|y|_\infty > B } H(y,x) \mu(y,x) dydx = \int \lt(\int_{|y|_\infty > B } H(y,x) \mu(y|x) dy \rt) \mu_X^\ast(x) dx  < \eps,
\end{equation}
where
\begin{equation}\nonumber
    B = 2\sigma_{\max} \lt( \log \frac{1}{\eps} + D \log \frac{\sigma_{\max}}{\sigma_{\min}} + \frac{K^2D}{2\sigma_{\max}^2} + \log 2D \rt)^{1/2}.
\end{equation}
For each $(i,j)$ define 
$$
l_{ij}(y,x) = \max\lt\{ p_{g_i, \sigma_j}(y,x) - \eps/2, 0 \rt\} \qquad\text{and}\qquad u_{ij}(y,x) = \min\lt\{ p_{g_i, \sigma_j}(y,x) + \eps/2, H(y,x)  \rt\}.
$$
It follows that
\begin{equation}
    \begin{split}
            &\int\int \lt\{ u_{ij}(y,x) - l_{ij}(y,x) \rt\}  \mu_X^\ast(x) dy dx 
            \\
            \le &\int \int_{|y|_\infty \le B} \eps  \mu_X^\ast(x) dy dx + \int \int_{|y|_\infty > B} H(y,x) \mu_X^\ast(x) dy dx 
            \\
            \le &\lt\{ (2B)^D + 1 \rt\} \eps. 
    \end{split}   
\end{equation}
Denote $\delta^2 := \lt\{ (2B)^D + 1 \rt\}$. With $d_H^2(u_{ij}, l_{ij}) \le d_1(u_{ij}, l_{ij}) $, we have
\begin{equation}\label{eq: t9}
    \calN_{[]}(\delta, \calP, d_H) \le \calN_{[]}(\delta^2, \calP, d_1) \le N_1N_2 \le \frac{\sigma_{\max} - \sigma_{\min}}{\eta_2} \calN(\eta_1, \F, \||\cdot|_\infty \|_\infty). 
\end{equation}

It is possible to write
\begin{equation}\nonumber
    \delta^2 = \eps \le C_1(\sigma_{\max}, D) \lt[ \eps (\log \eps^{-1})^{D/2} + \eps C_2(K)  + \eps \lt( \log \frac{\sigma_{\max}}{\sigma_{\min}} \rt)^{D/2}    \rt],
\end{equation}
where $C_1(\sigma_{\max}, D)$ and $C_2(K)$ is a constant. There exists small enough $\eps_*(D)$ such that for all $\eps \in (0,\eps_*]$ 
$$
\delta^2 \le C_3(\sigma_{\max}, D, K) \sqrt{\eps} \lt( \log \frac{\sigma_{\max}}{\sigma_{\min}} \rt)^{D/2}. 
$$
Consequently, there exists $\delta_* = \delta_*(D)$, such that for all $\delta \le \delta_*$, we have
\begin{equation}\nonumber
    C_3^2(\sigma_{\max}, K, D)\delta^4 \lt( \log \frac{\sigma_{\max}}{\sigma_{\min}} \rt)^{-D} \le \eps.
\end{equation}
It lead us to, for all $\delta \le \delta_*$
\begin{equation}\label{eq: t10}
    \eta_1 \ge \frac{ c_1^{-1}C_3^2 \sigma_{\min}^{D+3} \delta^4}{\sigma_{\min} \{\log(\sigma_{\max}/\sigma_{\min})\}^D} \ge c \sigma_{\min}^{D+3} \delta^4,
\end{equation}
where $c(\sigma_{\max}, K, D)$ is a constant. We use the fact that $\sigma_{\min} \{\log(\sigma_{\max}/\sigma_{\min})\}^D$ is bounded above by some constant depending only upon $\sigma_{\max}$ as $\sigma_{\min} \le 1$. Similar to \eqref{eq: t10}, it is possible to write for all $\delta > \delta_*$
\begin{equation}\label{eq: t11}
    \eta_2 \ge c'\sigma_{\min}^{D+2} \delta^4, \qquad\qquad \text{ for all } \delta \le  \delta_*,
\end{equation}
where $c'(\sigma_{\max}, K, D)$ is some constant.

\
\\
The result now follows directly \eqref{eq: t11} and \eqref{eq: t10} with \eqref{eq: t9}.
\end{proof}

\section{Proof of Theorem~\ref{thm: main}}\label{sec: proof thm: main}
\begin{proof}
    Choose four absolute constants $c_1, \ldots, c_4$ as in Theorem 1 of \citet{wongandshen}. Define $c$ and $C$ in the statement of \Cref{lemma: Entropy bound}. The proof closely follows \citet{linsparse}. We have therein the proof of Theorem 3 that
    \begin{equation} \label{eq: integrated brackett}
    \begin{split}
                &\int_{\eps^2/2^8}^{\sqrt{2}\eps} \sqrt{\log \calN_{[]} (\delta/c_3, \p, d_H) }  d\delta 
                \\
                \le & \, \sqrt{2} \eps \sqrt{\xi A + (D+3)(s+1) \log \sigma_{min}^{-1} + c_5 \xi } + \sqrt{2} \eps \sqrt{4(\xi+1)} \sqrt{\log (2^8/\eps^2)},
    \end{split}
    \end{equation}
    for every $\eps \le \sqrt{2} \le c_3 \delta_*/\sqrt{2}$, where $c_5 = c_5(c, C, c_3)$. Observe that $c_4 \sqrt{n} \eps_n^2$ is upper bound to \eqref{eq: integrated brackett} and Eq. (3.1) of \citet{wongandshen} is satisfied.

    \
    \\
    Using B.12 of \citet{ghosalvandervaart2017}, we have
\begin{align*}
    K(p_{G_*, \sigma_*}, p_{g, \sigma_*}) &\le \int \int K\Big( N\lt(G_*(z,x), \sigma_*^2\rt) ,N\lt(g(z,x), \sigma_*^2  \rt) \Big)  \mu_X^\ast(x) \thinspace dx \thinspace dP_Z(z) 
    \\
    &= \int \int \frac{|G_*(z,x) - g(z,x)|^2_2}{2\sigma_*^2}  \mu_X^\ast(x) \thinspace dx \thinspace dP_Z(z) \le \frac{D \delta_{\app}^2}{2 \sigma_*^2} =: \delta_n.
\end{align*}
One may easily see that
$$
\int \lt( \log \frac{\phi_\sigma(x)}{\phi_\sigma(x-y)} \rt)^2 \phi_\sigma(x) dx = \int \frac{|y|^4_2 + 4|x^Ty|^2}{4\sigma^2} \phi_\sigma(x) dx \le \frac{|y|_2^4}{4\sigma^2} + |y|_2^2 \int \frac{|x|_2^2}{\sigma^2} \phi_\sigma(x) dx.
$$
Combining this with Example B.12, (B.17) and Exercise B.8 of \citet{ghosalvandervaart2017}, we have
\begin{align*}
    &\int \int \lt( \log \frac{p_{G_*, \sigma_*}(y|x)}{p_{g,\sigma_*}(y|x)} \rt)^2 \, dP_*(y|x) \, \mu_X^\ast(x) dx
    \\  
     \le &\int \int \int \lt( \log \frac{\phi_\sigma(y - G_*(z,x)}{\phi_\sigma(y - G(z,x)} \rt)^2 \phi_\sigma(y - G_*(z,x)) \,  dy  \thinspace dP_Z(z) \thinspace \mu_X^\ast(x) dx
    \\
    \le &\frac{D^2 \delta^4_{\app}}{4\sigma_*^2} + D\delta^2_{\app} \int \frac{|x|^2_2}{\sigma_*^2} \phi_{\sigma_*} (y) dy + \frac{2D\delta^2_{\app}}{\sigma_*^2} \le c_7\frac{\delta_{\app}^2}{\sigma_*^2} =: \tau_n,
\end{align*}
where $c_7 = c_7(D)$. We are using $\delta_n$ and $\tau_n$, although they are independent of $n$, for notational consistency with Theorem 4 of \citet{wongandshen}. Let $\eps_n^* = \eps_n \vee \sqrt{12 \delta_n}$. Then, using Theorem 4 of \citet{wongandshen}, we have
$$
P_* \lt( d_H(\widehat{p}, p_*) > \eps_n \rt) \le 5 e^{-c_2 n \eps_n^{*2}} + \frac{\tau_n}{n \delta_n} = 5 e^{-c_2 n \eps_n^{*2}} + \frac{2c_7^2}{Dn}.
$$
The proof is complete after redefining constants.
\end{proof}

\section{Proofs of Corollary~\ref{cor: both case}}\label{sec: proof cor: both case}
\begin{proof}
    For the sparse case in \ref{cor: both case}.1, utilizing the entropy bound from \eqref{eq: entropy Hieber}, we observe that
\[ 
\xi\{A + \log(n/\sigma_{\min})\} \asymp \delta_{\app}^{-t_*/\beta_*} \log^3(\delta^{-1}_{\app}),
\]
which naturally leads to the required convergence rate. 

    Similarly for the fully connected case \ref{cor: both case}.2, utilizing the entropy bound from \eqref{eq: entropy Lin} , we observe that
        \[ 
        \xi\{A + \log(n/\sigma_{\min})\} \asymp \delta_{\app}^{-t_*/\beta_*} \log^3(\delta^{-1}_{\app}),
        \]
    which naturally leads to the required convergence rate.
\end{proof}

\section{Proof of Theorem~\ref{thm: wasserstein general}}\label{sec: proof thm: wasserstein general}

\begin{proof}
    It is suffice to assume that $\eps$ and $\sigma_* \sqrt{\log \eps^{-1}}$ are sufficiently small. If not, let $\eps + \sigma_* \sqrt{\log \eps^{-1}} \ge c_0$, where $c_0(K,D, \mathsf{r}_\ast)$. Then \Cref{thm: wasserstein general} holds trivially by taking a large enough constant depending just on $D$, $K$, and $\mathsf{r}_\ast$.

    Let $V \sim Q(\cdot|X=x)$, $V_* \sim Q(\cdot|X=x)$, $\beps \sim \normal(0_D, \sigma^2 \mathbb{I}_d) $ and $\beps_* \sim \normal(0_D, \sigma_*^2 \mathbb{I}_d) $ be independent with underlying probability density $\nu$. We truncate the random variable $\beps$ and $\beps_*$ componentwise as $\lt(\beps_K\rt)_j = \max\{-K, \min\{K, \beps_j\} \}$ and $\lt(\beps_{*K}\rt)_j = \max\{-K, \min\{K, \lt(\beps_{*}\rt)_j\} \}$ respectively. We denote $P_{g,\sigma}$ as $P$, $Q_g$ as $Q$, $\wtilde{P}$ as distribution of $V + \beps_K$ and $\wtilde{P}_*$ as the distribution of $V_* + \beps_{*K}$. One may note that $W_1(\wtilde{P}_*, Q_*) \le W_2(\wtilde{P}_*, Q_*) \le \sqrt{\E\big[ |\beps_{*K}|_2^2 \big]} \le \sqrt{\E\big[ |\beps_{*}|_2^2 \big]} \le \sigma_* \sqrt{D}$. Similarly, $W_1(\wtilde{P}, Q) \le \sigma \sqrt{D}$. The $\ell_1$ diameter of $[-2K, 2K]^D$, where the support of $\wtilde{P}$ and $\wtilde{P}_*$, is 4KD. Observe that
    $$
    W_1\lt(\wtilde{P}_*,\wtilde{P}\rt) \le 4\,KD\, d_1\lt(\wtilde{P}_*,\wtilde{P}\rt) \le 4\,KD\, d_1(P_*,P) \le 8\,KD\, d_H(P_*,P), 
    $$

where the first inequality follows from Theorem 4 of \citet{gibbs2002}, the second inequality follows from the fact the distance between two truncated distributions is always lesser than the original distributions and the last inequality follows from $d_1 \le 2d_H$. Hence,
$$
W_1\lt(Q_*, Q \rt) \le W_2\lt( Q_*, \wtilde{P}_* \rt) +  W_1\lt(\wtilde{P}_*,\wtilde{P}\rt) + W_2\lt(\wtilde{P}, Q \rt) \le \sigma_* \sqrt{D} + 8\, KD\, \eps + \sigma \sqrt{D}.
$$
Now it is suffice to show that $\sigma \le c \;\sigma_* \sqrt{\log \eps^{-1}}$, where $c= c(D, K, r*)$ is a constant, because we have assumed that $\eps$ is small enough. We establish this in the rest of the proof. Let $t_* = \lt[ 2 \,\sigma_*^2\, D \log \lt(\frac{2D}{\eps}\rt) \rt]^{1/2}$. Observe that
$$
\int_{|x|_2 >t_*} \phi_{\sigma_*}(x) dx \le \int_{|x|_\infty > t_*/\sqrt{D}} \phi_{\sigma_*}(x) dx \le 2\,D e^{-t_*^2/2D\sigma^2} \le \eps.
$$
Let $\M_*^{t_*} = \M_* \oplus \calB_{t_*}(0_D)$. We may write
\begin{equation}\label{eq: manifold boundary prob}
    \begin{aligned}
        &1 - P_*\lt(\M_*^{t_*}\rt) = \nu\lt( Y_* + \beps_* \notin \M_*^{t_*}  \rt) \le \nu\lt(|\beps_*|_2>t_*  \rt) 
        \\
        \implies& P\lt(\M_*^{t_*}\rt) \ge 1 - 2\eps,
    \end{aligned}
\end{equation}
the implication in the last line follows from $\sup_{B}\lt| P(B) - P_*(B) \rt| \le  d_H(P, P_*) \le \eps$. For the sake of contradiction, let $\sigma \in \lt[2t_*, r^*/2\rt]\cup \lt( \mathsf{r}_\ast/2, \infty\rt)$ ($t_*$ is sufficiently small, from the assumption we made at the beginning of this proof). If $\sigma > \mathsf{r}_\ast/2$, then
$$
2\eps \ge 1 - P\lt(\M_*^{t_*}\rt) \ge 1 - P\lt( [-K , K]^D \rt) \ge c_2(K, D, r*)
$$
where $c_2$ is some positive constant. It is a contradiction following from the smallness of $\eps$. Lets make a claim that if  $\sigma \in \lt[2t_*, \mathsf{r}_\ast/2\rt]$, then for every $y \in \R^D$, there is some $z \in \R^D$ such that $|z-y|_2 \le \sigma$ and $\calB_{\sigma/2}(z) \cap \M_*^{t_*} = \emptyset$.

\
\\
Following from the claim, we have
$$
\nu\lt(  Y + \beps \notin \M_*^{t_*} \big| Y = y\rt) \ge \nu\lt( \beps \in \calB_{\sigma/2}(z - y)  \rt).
$$
Since $|z-y|_2 \le \sigma$, the right hand side is bounded below by a positive constant depending just on $D$ which is again a contradiction of \eqref{eq: manifold boundary prob}. This proves the assertion made in the theorem.

\
\\
The proof of the claim is divided into three cases. Let $\rho\lt( y, \M_*\rt) = \inf\{|y - y'|_2: y'\in \M_*\} $ be the $\ell_2$ set distance.

{\bf Case 1.} $\rho(y, \M_*) \ge \sigma:$ We may choose $z=y$.

{\bf Case 2.} $\rho(y, \M_*) \in (0, \sigma):$ Let $y_0$ be the unique Euclidean projection of $y$ onto $\M_*$. Such a unique projection exists because $\sigma < \mathsf{r}_\ast$ is within the reach and $y\in \M_*$, since $\M_*$ is closed. Suppose $y_t = y_0 + t (y - y_0)$. We shall define two continuous functions $d_0(t) = |y_t - y_0|_2$ and $d(t) = \rho(y_t, \M_*)$. It is obvious that $d(t) \le d_0(t)$. For $t\in \big[0,1+\sigma/|y-y_0|_2\big]$, $d_0(t) \le d(t)$ because $y_0$ is the unique projection for all the points that lie on the line segment including the farthest point with $t = 1+\sigma/|y-y_0|_2$. Otherwise, say $d(t) = \rho(y_t, z)$ and 
$$
|y - y_0|_2 = |y - y_t|_2 + |y_t - y_0|_2 > |y - y_t| + |y_t - z| \ge |y - z|_2 
$$
which contradicts $y_0$ being a unique projection. The claim holds for the point $z = y_{1+\sigma/|y-y_0|_2}$. To see this, observe $|z-y| = \sigma$ and $\calB_{\sigma/2}(z) \cap \M_*^{t_*} = \emptyset$ because $t_* \le \sigma/2$ and the ball $\calB_{\sigma/2}(z) \subset \M_*^{\mathsf{r}_\ast}$ is within the reach of the manifold.

{\bf Case 3.} $\rho(y, \M_*) = 0:$ Because $\M_*$ has empty interior, for all $\gamma>0$, we always find a point $y_\gamma$, which in $\calB_\gamma(y)$ which away from $\M_*$. For small enough $\gamma$, we reduce to case 2 by taking $\gamma \to 0$, the limit point of $y_\gamma$ has the required behavior. 

\end{proof}
    
\section{Proof of Corollary~\ref{cor: perturbation}}\label{sec: proof cor: perturbation}
\begin{proof}
The effective noise variance after the perturbation would be
\begin{align*}
    \wtilde{\sigma}_* = n^{-\alpha} + n^{-\beta_*/2(\beta_* + t_*)} \asymp \begin{cases}
        n^{-\alpha}, \qquad &\alpha < \beta_*/\{2(\beta_* + t_*)\} \\
        n^{\beta_*/2(\beta_* + t_*)}, \qquad &\text{otherwise}.
    \end{cases}
\end{align*}
Following this and the \Cref{thm: wasserstein general}, for the rate we have
\begin{align*}
    \eps^*_n + \sigma_* \sqrt{\log ((\eps^*_n)^{-1})} &\asymp  \lt(n^{-\frac{\beta_* - t_*\alpha}{2\beta_* + t_*}} + n^{-\alpha} \rt)\log^{2} (n) 
    \\
    &\asymp 
    \begin{cases}
        n^{-\frac{\beta_* - t_*\alpha}{2\beta_* + t_*}}\, \log^{2} (n) , \quad &\text{if}\; \alpha < \beta_*/\{2(\beta_* + t_*)\} ,    
        \\
        n^{-\frac{\beta_*}{2(\beta_* + t_*)}}\, \log^{2} (n) , \quad &\text{otherwise}.
    \end{cases}
\end{align*}

\end{proof}

\section{Proof of Theorem~\ref{thm: disparity}}\label{sec: proof thm: disparity}
\begin{proof}
 With $m = \lceil \log_2 (n) \rceil$ and $N = \lt(n^{(\beta_Z^{-1}d + \beta_X^{-1}p)\lt[1+\alpha(\beta_Z^{-1}d + \beta_X^{-1}p)\rt]/\lt[2 + \beta_Z^{-1}d + \beta_X^{-1}p \rt] }\rt)$  in \Cref{thm.approx_network_One_fct2}, we can find a network $G$ with the mentioned architecture such that
$$
\||G - G_* |_\infty \|_\infty \le \delta_{\app}.
$$
Following the entropy bound from \eqref{eq: entropy Hieber}, we have
\begin{align*}
      \log \calN(\delta, \F_s, \||\cdot|_\infty\|_\infty ) &\lesssim sL\, \{ \log(rL) + \log \delta^{-1} \} 
      \\
      &\lesssim \delta^{-(\beta_Z^{-1}d + \beta_X^{-1}p )}_{\app} \log^2 \delta^{-1}_{\app} \lt\{\ \log \lt( \delta^{-1}_{\app} \log \lt( \delta^{-1}_{\app}  \rt) \rt)  + \log \lt( \delta^{-1}_{\app}  \rt) \rt\}. 
\end{align*}
The rest directly follows from the \Cref{thm: main}
\end{proof}

\section{Approximation properties of the sparse and fully connected DNNs}\label{sec: approx old}
The approximability of the sparse network is detailed in \Cref{lemma: DNN approximation}.1, which restates Lemma 5 from \citet{linsparse}. For the fully connected network, \Cref{lemma: DNN approximation}.2 demonstrates its approximation capabilities, derived directly from Theorem 2 and the proof of Theorem 1 in \citet{kohlerlanger}. Additionally, the inclusion of the class $\G$ in the fully connected setup is supported by the discussion in Section 1 of \citet{kohler2020discussion}.
\begin{lemma}\label{lemma: DNN approximation}
    Suppose that $G_* \in \G$. Then, for every small enough $\delta \in (0,1)$,
    
    \begin{enumerate}
        \item there exists a sparse network $G \in \F_s = \F_s\lt(L, r, s, K \vee 1\rt)$ with $L \lesssim \log \delta^{-1}$, $r \lesssim \delta^{-t_*/\beta_*}$, $s \lesssim \delta^{-t_*/\beta_*} \log \delta^{-1} $ satisfying $\| |G - G_*|_\infty\|_\infty \le \delta$.
        \item there exists a fully connected network $G \in \F_c$ with $L \lesssim \log \delta^{-1}$, $r \lesssim \delta^{-t_*/2\beta_*}$, $B \lesssim \delta^{-1}$ satisfying $ \lt\| \lt|G - G_* \rt|_\infty \rt\|_\infty \le \delta$.
    \end{enumerate}
\end{lemma}

\section{A new approximation result for functions with smoothness disparity}\label{sec: approx}

In this section, we prove the approximability of the sparse neural network for the H\"older class of function $f\in {\mathcal{H}}_{{r}, {r_\prime}}^{\beta, \beta_\prime} \left(D,D_\prime, K \right)$.

\begin{theorem}\label{thm.approx_network_One_fct2}
Let $f\in {\mathcal{H}}_{{r}, {r_\prime}}^{\beta, \beta_\prime} \left( [0,1]^{{r}}, [0,1]^{{r_\prime}}, K \right)$. Denote $\rsum = r + r_\prime$ and $\bsum = \beta + \beta_\prime$. Then for any integers $m \geq 1$ and $N \geq  (\bsum + 1)^{\rsum} \vee (K+1)e^{\rsum},$ there exists a network
$$
\widetilde f \in \mF_s\big(L, \big(\rsum, 6(\rsum +\lceil \bsum \rceil)N, \ldots, 6(\rsum +\lceil \bsum \rceil)N ,1 \big), s, \infty\big)$$
with depth 
$$
L=8+(m+5)\left(1+ \left\lceil \log_2 \big(\rsum \vee \bsum \big) \right\rceil \right)
$$ 
and the number of parameters
\begin{align*}
	s\leq 109 \big(\rsum + \bsum + 1 \big)^{3+ \rsum} N (m+6),
\end{align*}
such that 
\begin{align*}
	\| \widetilde f - f\|_{L^\infty([0,1]^{\rsum})}\leq  (2K+1)\left( 1+ {\mathsf r}_{\mathrm{sum}}^2 + \beta_{\mathrm{sum}}^2 \right) 6^{\rsum}\, N\, 2^{-m}+ K\, 3^{\rsum/(\beta^{-1}r + \beta_\prime^{-1}r_\prime)}\, N^{-1/(\beta^{-1}r + \beta_\prime^{-1}r_\prime)}.
\end{align*}
\end{theorem}
We denote $\wtilde{\beta} = (\beta + \beta_\prime)^{-1} \beta \beta_\prime$ and $\wtilde{r} = (\beta + \beta_\prime)^{-1}(r \beta + r_\prime \beta_\prime)$.
Before presenting the proof of \Cref{thm.approx_network_One_fct2}, we formulate some required results. 

\
\\
We follow the classical idea of function approximation by local Taylor approximations that have previously been used for network approximations in \cite{yarotski2017} and \cite{schmidtheiber}. For a vector $\ba \in [0,1]^r$ define 
\begin{align}
	P_{\ba, \bb}^{\beta, \beta_\prime} f(\bu, \bv) = \sum_{ \substack{ 0 \leq |\balpha|< \beta \\ 0 \leq |\balpha_\prime|< \beta_\prime }} (\partial^{\balpha + \balpha_\prime} f) (\ba, \bb) \frac{(\bu - \ba)^{\balpha} (\bv - \bb)^{\balpha_\prime}}{\balpha !\, \balpha_\prime ! }.
	\label{eq.Taylor_approx}
\end{align}
We use the notation the $\bu = (u^{(j)})_j$ to represent the component of the vector when the index $j$ is well understood. Accordingly we have $\bv = (v^{(j)})_j$, $\ba = (a^{(j)})_j$ and $\bb = (b^{(j)})_j$ . By Taylor's theorem for multivariate functions, we have for a suitable $\xi \in [0,1],$ 
\begin{align*}
	f(\bu,\bv)= &\sum_{ \substack{\balpha : |\balpha| < \beta-1 \\ \balpha_\prime : |\balpha_\prime| < \beta_\prime-1}} (\partial^{\balpha + \balpha_\prime} f) (\ba, \bb) \frac{(\bu - \ba)^{\balpha} (\bv - \bb)^{\balpha_\prime}}{\balpha !\, \balpha_\prime ! }
    \\
	\, + 
	& \sum_{ \substack{\beta - 1 \leq |\balpha | < \beta \\ \beta_\prime - 1 \leq |\balpha_\prime | < \beta_\prime}  } (\partial^{\balpha + \balpha_\prime} f) (\ba + \xi (\bu-\ba) , \bb + \xi (\bv-\bb)) \frac{(\bu - \ba)^{\balpha} (\bv - \bb)^{\balpha_\prime}}{\balpha !\, \balpha_\prime ! }.
\end{align*}
We have $|(\bu - \ba)^{\balpha}| =\prod_{j=1}^r |u_j - a_j |^{\alpha^{(j)}} \leq |\bu-\ba|_\infty^{|\balpha|}$ and $|(\bv - \bb)^{\balpha_\prime}| =\prod_{j=1}^{r_\prime} |v_j - b_j |^{\alpha_\prime^{(j)}} \leq |\bv-\bb|_\infty^{|\balpha_\prime|}$. Consequently, for $f\in {\mathcal{H}}_{{r}, {r_\prime}}^{\beta, \beta_\prime} \left( [0,1]^{r}, [0,1]^{r_\prime}, K \right)$,
\begin{align}
	&\big | f(\bu, \bv) - P_{\ba, \bb}^{\beta, \beta_\prime} f(\bu, \bv)  \big| \notag 
    \\\label{eq.Hoeld_fct_approx} 
    \leq & \sum_{ \substack{\beta - 1 \leq |\balpha | < \beta \\ \beta_\prime - 1 \leq |\balpha_\prime | < \beta_\prime}  } \left( \partial^{\balpha + \balpha_\prime} f (\ba + \xi (\bu-\ba) , \bb + \xi (\bv-\bb) ) - \partial^{\balpha + \balpha_\prime}f(\ba, \bb) \right) \frac{(\bu - \ba)^{\balpha} (\bv - \bb)^{\balpha_\prime}}{\balpha !\, \balpha_\prime ! }
     \\
    \leq 
	&  K \left(|\bu-\ba|_\infty^{\beta} \vee |\bv-\bb|_\infty^{\beta_\prime} \right)  \notag
\end{align}
We may also write \eqref{eq.Taylor_approx} as a linear combination of monomials 
\begin{align}
	P_{\ba, \bb}^{\beta, \beta_\prime} f(\bu, \bv) = \sum_{ \substack{ 0 \leq |\bgamma| < \beta \\ 0 \leq |\bgamma_\prime| < \beta_\prime }}  c_{\bgamma, \bgamma_\prime} \bu^{\bgamma} \bv^{\bgamma_\prime} ,
	\label{eq.Pf_monomial_sum}
\end{align}
for suitable coefficients $c_{\bgamma, \bgamma_\prime}.$ For convenience, we omit the dependency on $\ba$ and $\bb$ in $c_{\bgamma, \bgamma_\prime}.$ Since 
$$
\partial^{\bgamma, \bgamma_\prime } P_{\ba, \bb}^{\beta, \beta_\prime} f(\bu, \bv) \, |_{(\bu=0, \bv = 0)} = \bgamma ! \, \bgamma_\prime ! \, c_{\bgamma, \bgamma_\prime},
$$ we must have
\begin{align*}
	c_{\bgamma, \bgamma_\prime} = \sum_{ \substack{ \bgamma \leq \balpha \&  |\balpha|< \beta \\ \bgamma_\prime \leq \balpha_\prime \&  |\balpha_\prime|< \beta_\prime }} (\partial^{\balpha  + \balpha_\prime} f) (\ba, \bb) \frac{(- \ba)^{\balpha -\bgamma} \, (- \bb)^{\balpha_\prime -\bgamma_\prime}}{\bgamma ! \, \bgamma_\prime ! \, (\balpha -\bgamma)! \, (\balpha_\prime -\bgamma_\prime)!}.
\end{align*} 
Notice that since $\ba \in [0,1]^{r}$, $\bb \in [0,1]^{r_\prime}$,  and $f \in {\mathcal{H}}_{{r}, {r_\prime}}^{\beta, \beta_\prime} \left( [0,1]^{r},[0,1]^{r_\prime}, K \right)$, 
\begin{align}
	|c_{\bgamma \bgamma_\prime}| \leq K/(\bgamma !\, \bgamma_\prime !) \quad \text{and} \ \ \sum_{ \substack{\bgamma \geq \bzero \\ \bgamma_\prime \geq \bzero }} |c_{\bgamma, \, \bgamma_\prime}| \leq K \prod_{i=1}^{r} \prod_{j=1}^{r_\prime} \sum_{\gamma^{(i)} \geq 0} \sum_{\gamma_\prime^{(j)} \geq 0}  \frac{1}{\gamma^{(i)}!} \,\frac{1}{\gamma_\prime^{(j)}!} = Ke^{r+r_\prime},
	\label{eq.cgamma_bd}
\end{align}
where $\gamma = (\gamma^{(1)}, \dots, \gamma^{(r)})$ and $\gamma_\prime = (\gamma_\prime^{(1)}, \dots, \gamma_\prime^{(r_\prime)})$.


\
\\
Consider the set of grid points 
\begin{align*}
    \bD(M) :=\{\bu_{\bell^\1} = & (\ell^\1_j/M_1)_{j=1,\ldots,r}  \text{ and } \bv_{\bell^\2} = (\ell^\2_j/M_2)_{j=1,\ldots,r_\prime} 
    \\
    : \, &\bell^\1 =(\ell^\1_1,\ldots, \ell^\1_r) \in \{0,1,\ldots, M_1 \}^{r},
    \\
    &\bell^\2 =(\ell^\2_1,\ldots, \ell^\2_r) \in \{0,1,\ldots, M_2 \}^{r_\prime}, M_1 = M^{\wtilde{\beta}/\beta}, M_2 = M^{\wtilde{\beta}/\beta_\prime} \}.    
\end{align*}

The cardinality of this set is $(M_1 + 1)^{r} \cdot (M_2 + 1)^{r_\prime} .$ We write $\bu_{\bell^\1} = (u^{(j)}_{\bell^\1})_{j=1, \ldots, r}$ and $\bv_{\bell^\2}=(v^{(j)}_{\bell^\2})_{j=1, \ldots, r_\prime}$ to denote the components of $\bu_{\bell^\1}$ and $\bv_{\bell^\2}$ respectively. With slight abuse of notation we denote $\bw = (\bu, \bv) = (u^\1, \ldots, u^{(r)}, v^\1, \ldots, v^{(r_\prime)})$, $\bell = (\bell^\1, \bell^\2) = (\ell^\1_1,\ldots, \ell^\1_{r}, \ell^\2_1, \ldots, \ell^\2_{r_\prime} )$ and $\bw_\bell = (w^{(j)}_\bell)_{j=1, \ldots, r + r_\prime}  = (\bu_{\bell^\1}, \bv_{\bell^\2}) = (u^{(1)}_{\bell^\1}, \ldots, u^{(r)}_{\bell^\1}, v^{(1)}_{\bell^\2}, \ldots, u^{(r_\prime)}_{\bell^\2}  )$.
Define
\begin{align*}
	&P^{\beta, \beta_\prime} f(\bu, \bv)
    \\=\,&P^{\beta, \beta_\prime} f(\bw) 
    \\
    := & \sum_{\bw_{\bell} \in \bD(M)} P_{\bw_{\bell}}^{\beta, \beta_\prime} f (\bw) \prod_{j=1}^{r+r_\prime} (1- M_j |w^{(j)} - w_{\bell}^{(j)}|)_+ 
    \\
    = &\sum_{\bu_{\bell^\1}, \bv_{\bell^\2}  \in \bD(M)} P_{ \bu_{\bell^\1}, \bv_{\bell^\2} }^{\beta, \beta_\prime} f (\bu, \bv) \left(\prod_{j=1}^{r} (1- M_1 |u^{(j)} - u_{\bell^\1}^{(j)}|)_+ \right) \left( \prod_{j=1}^{r_\prime} (1- M_2 |v^{(j)} - v_{\bell^\2}^{(j)}|)_+ \right),
\end{align*}
where $M_j = M_1$ for $j=1, \ldots, r$ and $M_j = M_2$ for $j=r + 1, \ldots, r + r_\prime$.

\medskip
\begin{lemma}\label{lem.Hoeld_approx}
If $f\in {\mathcal{H}}_{{r}, {r_\prime}}^{\beta, \beta_\prime} \left( [0,1]^{r}, [0,1]^{r_\prime}, K \right)$, then $\| P^{\beta, \beta_\prime} f - f \|_{L^\infty[0,1]^{r + r_\prime}} \leq K M^{-\wtilde{\beta}}.$
\end{lemma}

\begin{proof}
Since for all $\bw=(w^\1, \ldots, w^{(r+r_\prime)}) \in [0,1]^{r+r_\prime},$ 
\begin{align}
    \sum_{\bw_{\bell} \in \bD(M)} \prod_{j=1}^{r+r_\prime} (1- M_j |w^{(j)} - w_{\bell}^{(j)}|)_+ = \prod_{j=1}^{r + r_\prime} \sum_{\ell=0}^{M_j}(1- M_j |w^{(j)} - \ell/M_j|)_+ =1,
	\label{eq.kernel_sum_one}
\end{align}
we have 
\begin{align*}
    f(\bw) &= f(\bu, \bv) 
    \\
    &= \sum_{ \substack{\bu_{\bell^\1}, \bv_{\bell^\2} \in \bD(M):
\\
\|\bu- \bu_{\bell^\1}\|_\infty \leq 1/M_1 \\ \|\bv- \bv_{\bell^\2}\|_\infty \leq 1/M_2}} \, f(\bu, \bv) \left(\prod_{j=1}^{r} (1- M_1 |u^{(j)} - u_{\bell^\1}^{(j)}|)_+ \right) \left( \prod_{j=1}^{r_\prime} (1- M_2 |v^{(j)} - v_{\bell^\2}^{(j)}|)_+ \right)    
\end{align*}

and with \eqref{eq.Hoeld_fct_approx},
\begin{align*}
	\big | P^{\beta, \beta_\prime} f (\bu, \bv) - f(\bu, \bv) \big| 
	&\leq \max_{\substack{\bu_{\bell^\1}, \bv_{\bell^\2} \in \bD(M): \\ \|\bu- \bu_{\bell^\1}\|_\infty \leq 1/M_1 \\ \|\bv- \bv_{\bell^\2}\|_\infty \leq 1/M_2}} \big| P_{ \bu_{\bell^\1}, \bv_{\bell^\2} }^{\beta, \beta_\prime} f (\bu, \bv) -  f(\bu, \bv) \big| 
	\\
    &\leq K \left(M_1^{-\beta} \vee M_2^{-\beta_\prime} \right) = K\,M^{-\wtilde{\beta}}.
\end{align*}
\end{proof}

\
\\
In the next few steps, we describe how to build a network that approximates $ P^{\beta, \beta_\prime} f $.

\medskip
\begin{lemma} \label{lem.hat_fct_mult}
Let $M,m,$ be any positive integer. Denote $M_1 = M^{\wtilde{\beta}/\beta}$, $M_2 = M^{\wtilde{\beta}/\beta_\prime}$, $\mathsf{M} = (M_1 + 1)^{r}(M_2 + 1)^{r_\prime}$ and $\rsum =  r + r_\prime$. Then there exists a network 
$$
\Hat^{\rsum} \in \mF\left( 2 + (m+5) \lceil \log_2 (\rsum) \rceil , \rsum, 2\rsum\mathsf{M}, \rsum \mathsf{M}, 6\rsum \mathsf{M},  \ldots, 6 \rsum \mathsf{M}, \mathsf{M}), s, 1\right)
$$
with $s\leq 37 \rsum^2 \mathsf{M} (m+5) \lceil\log_2 (\rsum) \rceil,$ such that $\Hat^r \in [0,1]^{\mathsf{M}}$ and for any $\bu=(u^\1, \ldots, u^{(j)}) \in [0,1]^{r}$ and for any $\bv=(v^\1, \ldots, v^{(j)}) \in [0,1]^{r_\prime}$
\begin{align*}
	\Bigg| \Hat^{\rsum}(\bu,\bv) - \Bigg\{&\Big(\prod_{j=1}^{r} (1/M_1 -  |u^{(j)} - u^{(j)}_{\bell^\1}|)_+ \Big) \times
    \\
    &\Big(\prod_{j=1}^{r_\prime} (1/M_2 -  |v^{(j)} - v^{(j)}_{\bell^\2}|)_+\Big) \Bigg\}_{\bu_{\bell^\1}, \bv_{\bell^\2}  \in \bD(M)} \Bigg|_\infty
	\leq \rsum^2 2^{-m}.
\end{align*}
For any $\bu_{\bell^\1}, \bv_{\bell^\2}  \in \bD(M),$ the support of the function $(\bu, \bv) \mapsto (\Hat^{r + r_\prime}(\bu,\bv))_{\bu_{\ell^\1}, \bv_{\ell^\2}}$ is moreover contained in the support of the function
$$
(\bu, \bv) \mapsto \left\{\bigg(\prod_{j=1}^{r} (1/M -  |u^{(j)} - u^{(j)}_{\bell^\1}|)_+\bigg)\bigg( \prod_{j=1}^{r_\prime} (1/M -  |v^{(j)} - v^{(j)}_{\bell^\2}|)_+ \bigg) \right\}.
$$
\end{lemma}

\begin{proof}

{\bf Step~$1$:} (For $r + r_\prime =1$) Without loss of generality we consider the case when $r = 1$ and $r_\prime = 0$. 
We compute the functions $\{ (u^{(j)} - \ell/M_1)_+\}_{j=1, \ell = 0}^{r, M_1}$ and $\{ (\ell/M_1 - u^{(j)})_+\}_{j=1, \ell = 0}^{r, M_1}$ for the first hidden layer of the network. This requires $2r(M_1 +1)$ units (nodes) and $2r (M_1 + 1)$ non-zero parameters. 

For the second hidden layer we compute the functions $(1/M_1 -  |u^{(j)}-\ell/M_1|)_+ = (1/M_1 -  (u^{(j)} - \ell/M_1)_+ - (\ell/M_1 - u^{(j)})_+)_+$ using the output $(u^{(j)} - \ell/M_1)_+$ and $(\ell/M_1 - u^{(j)})_+$ from the output of the first hidden layer. This requires $r(M_1 + 1) + r_\prime(M_2 + 1)$ units (nodes) and $2 r (M_1 + 1)$ non-zero parameters. This proves the result for the base case when $r + r_\prime = 1$.

\
\\
{\bf Step~$2$:}
For $r + r_\prime >1,$ we compose the obtained network with networks that approximately compute the following
$$
\left\{\left(\prod_{j=1}^{r} (1/M_1 - |u^{(j)} - u^{(j)}_{\bell^\1}|)_+\right) \left(\prod_{j=1}^{r_\prime} (1/M_2 - |v^{(j)} - v^{(j)}_{\bell^\2}|)_+\right)\right\}_{\bu_{\bell^\1}, \bv_{\bell^\2}  \in \bD(M)}.
$$
For fixed $\bu_{\bell^\1}$ and $\bv_{\bell^\2}$, and from the use of \Cref{lem.multi_mult} there exist $\Mult_m^{r+r_\prime}$ networks in the class 
$$
\mF\left( 2 + (m+5) \lceil \log_2 (r + r_\prime) \rceil , (r + r_\prime, 2(r + r_\prime), r + r_\prime, 6(r+r_\prime), 6(r+r_\prime), \ldots, 6 (r + r_\prime), 1)\right)
$$ computing $(\prod_{j=1}^{r} (1/M_1 - |u^{(j)} - u_{\bell^\1}|)_+) \times (\prod_{j=1}^{r_\prime} (1/M_2 - |v^{(j)} - v_{\bell^\2}|)_+)$  up to an error that is bounded by $(r + r_\prime)^2\, 2^{-m}$. Observe that we have two extra hidden layers to compute $(1/M_1 - |u^{(j)} - u_{\bell^\1}|)_+)$ and $(1/M_2 - |v^{(j)} - v_{\bell^\2}|)_+)$ for fixed $\bu_{\bell^\1}$ and $\bv_{\bell^\2}$ respectively, before we enter into the multinomial computation by regime invoking \Cref{lem.multi_mult}. Observe that the number of parameters in this network is upper bounded by $37 (r + r_\prime)^2 (m+5) \lceil\log_2 (r + r_\prime) \rceil $. 

Now we use the \textit{parallelization} technique to have $(M_1 + 1)^{r}\cdot (M_1 + 1)^{r}$ parallel architecture for all elements of $\bD(M)$. This provides the existence of the network with the number of non-zero parameters bounded by $37 (r + r_\prime)^2 (M_1 + 1)^{r} (M_2 + 1)^{r_\prime} (m+5) \lceil\log_2 (r + r_\prime) \rceil $

By Lemma \ref{lem.multi_mult}, for any $\bx \in \real^r$, $\Mult_m^r(\bx)=0$ if one of the components of $\bx$ is zero. This shows that for any $\bu_{\ell^\1}, \bv_{\ell^\2} \in \bD(M)$, the support of the function $(\bu,\bv) \mapsto (\Hat^{r + r_\prime}(\bu,\bv))_{\bu_{\ell^\1}, \bv_{\ell^\2}}$ is contained in the support of the function $(\bu, \bv) \mapsto \left(\prod_{j=1}^{r} (1/M -  |u^{(j)} - u^{(j)}_{\bell^\1}|)_+ \prod_{j=1}^{r_\prime} (1/M -  |v^{(j)} - v^{(j)}_{\bell^\2}|)_+ \right)$.
\end{proof}
\begin{proof}[Proof of \Cref{thm.approx_network_One_fct2}]
All the constructed networks in this proof are of the form $\mF(L, \bp,s) =\mF(L, \bp,s, \infty)$ with $F=\infty$. Denote $M_1 = M^{\wtilde{\beta}/\beta}$, $M_2 = M^{\wtilde{\beta}/\beta_\prime}$, $\bsum = \beta + \beta_\prime$, and $\rsum = r + r_\prime$. Let $M$ be the largest integer such that $\mathsf{M} = (M_1+1)^{r} (M_2+1)^{r_\prime} \leq N$ and define $L^*:=(m+5)\lceil \log_2 (\bsum \vee \rsum) \rceil$.  Thanks to \eqref{eq.cgamma_bd}, \eqref{eq.Pf_monomial_sum} and Lemma \ref{lem.monomials}, we can add one hidden layer to the network $\Mon_{m, \bsum}^\rsum$ to obtain a network 
$$
Q_1 \in \mF\big( 2+L^* , (r, 6\lceil \beta\rceil C_{\rsum,\bsum}, \ldots, 6\lceil \beta\rceil C_{\rsum,\bsum}, C_{\rsum,\bsum}, \mathsf{M})\big),
$$
such that $Q_1(\bu, \bv) \in [0,1]^{\mathsf{M}}$ and for any $\bu\in [0,1]^{r}$ and for any $\bv \in [0,1]^{r_\prime}$
\begin{align}
	\Big | Q_1(\bu, \bv) - \Big( \frac{ P^{\beta, \beta_\prime} f (\bu, \bv) }{B}  + \frac 12\Big)_{\bu_{\bell^\1}, \bv_{\bell^\2}  \in \bD(M) } \Big|_\infty \leq \beta^2_{\mathrm{sum}} 2^{-m}
	\label{eq.Q1_bd}
\end{align}
with $B:= \lceil 2Ke^\rsum \rceil$. The total number of non-zero parameters in the $Q_1$ network is $6 \rsum (\bsum + 1) C_{\rsum,\bsum} + 42(\bsum + 1)^2C_{\rsum,\bsum}^2(L^*+1) +C_{\rsum,\bsum}\mathsf{M}$.

Recall that the network $\Hat^\rsum$ computes the products of hat functions (splines) $(\prod_{j=1}^{r} (1/M_1 - |u^{(j)} - u_{\bell^\1}|)_+) (\prod_{j=1}^{r_\prime} (1/M_2 - |v^{(j)} - v_{\bell^\2}|)_+)$ up to an error that is bounded by $\mathsf{r}^2_\mathrm{sum} 2^{-m}.$ It requires at most $37 \mathsf{r}^2_\mathrm{sum} N L^\ast$ active parameters. 
Observe that $C_{\rsum,\bsum} \leq (\bsum+1)^\rsum \leq N$ by the definition of $C_{r,\beta}$ and the assumptions on $N.$ By Lemma \ref{lem.hat_fct_mult}, the networks $Q_1$ and $\Hat^\rsum$ can be embedded into a joint parallel network $(Q_1, \Hat^\rsum)$ with $2+L^*$ hidden layers of size $(\rsum, 6(\rsum+\lceil \bsum\rceil)N,\ldots, 6(\rsum+\lceil \bsum\rceil) N, 2\mathsf{M})$. Using $C_{r,\beta} \vee (M+1)^r \leq N$ again, the number of non-zero parameters in the combined network $(Q_1, \Hat^r)$ is bounded by
\begin{align}
	\begin{split}
    & 6 \rsum (\bsum + 1) C_{\rsum,\bsum} + 42(\bsum + 1)^2C_{\rsum,\bsum}^2(L^*+1) +C_{\rsum,\bsum}\mathsf{M} + 37 \mathsf{r}^2_\mathrm{sum} N L^\ast \\ 
	&\leq 42 (\rsum + \bsum + 1)^2 C_{\rsum,\bsum} N (1+L^*)\\
	&\leq 84 (\rsum + \bsum + 1)^{3+\rsum} N (m+5),
	\end{split}
	\label{eq.para_Q1_bd}
\end{align}
where for the last inequality, we used $C_{\rsum,\bsum} \leq (\bsum+1)^\rsum,$ the definition of $L^*$ and that for any $x\geq 1,$ $1+\lceil \log_2(x)\rceil\leq 2+\log_2(x)\leq 2(1+\log(x)) \leq 2x.$

Next, we pair the $(\bu_{\bell^\1}, \bv_{\bell^\2})$-th entry of the output of $Q_1$ and $\Hat^r$ and apply to each of the $\mathsf{M}$ pairs the $\Mult_m$ network described in Lemma \ref{lem.mult}. In the last layer, we add all entries. By Lemma \ref{lem.mult} this requires at most $24(m+5)\mathsf{M} + \mathsf{M}\leq 25 (m+5)N$ active parameters for the $\mathsf{M}$ multiplications and the sum. Using Lemma \ref{lem.mult}, Lemma \ref{lem.hat_fct_mult}, \eqref{eq.Q1_bd} and triangle inequality, there exists a network $Q_2 \in \mF(2 + L^* + m+6, (\rsum, 6(\rsum + \lceil \bsum \rceil) N, \ldots, 6(\rsum +\lceil \bsum \rceil)N,1))$ such that for any $\bu\in [0,1]^{r}$ and for any $\bv \in [0,1]^{r_\prime}$

\begin{flalign*}
    	\Bigg | Q_2(\bu, \bv) - \sum_{ \bu_{\bell^\1}, \bv_{\bell^\2}  \in \bD(M) } \Big( \frac{ P^{\beta, \beta_\prime} f (\bu, \bv) }{B}  + \frac 12\Big) &\Big(\prod_{j=1}^{r} (1/M_1 - |u^{(j)} - u^{(j)}_{\bell^\1}|)_+\Big)& 
        \\
        &\Big(\prod_{j=1}^{r_\prime} (1/M_2 - |v^{(j)} - v^{(j)}_{\bell^\2}|)_+\Big) \Bigg|&
\end{flalign*}
\begin{flalign}\nonumber
	&\leq \sum_{\substack{\bu_{\bell^\1}, \bv_{\bell^\2} \in \bD(M): \\ \|\bu- \bu_{\bell^\1}\|_\infty \leq 1/M_1 \\ \|\bv- \bv_{\bell^\2}\|_\infty \leq 1/M_2}} (1+ \mathsf{r}^2_{\mathrm{sum}} + \beta^2_{\mathrm{sum}}) 2^{-m}&
    \\\label{eq.Q2_estimate}
	&\leq (1+ \mathsf{r}^2_{\mathrm{sum}} + \beta^2_{\mathrm{sum}}) 2^{r-m}.&
\end{flalign}
Here, the first inequality follows from the fact that the support of $(\Hat^{r + r_\prime}(\bu,\bv))_{\bu_{\ell^\1}, \bv_{\ell^\2}}$ is contained in the support of $\left(\prod_{j=1}^{r} (1/M -  |u^{(j)} - u^{(j)}_{\bell^\1}|)_+ \prod_{j=1}^{r_\prime} (1/M -  |v^{(j)} - v^{(j)}_{\bell^\2}|)_+ \right)$ (see Lemma \ref{lem.hat_fct_mult}). Because of \eqref{eq.para_Q1_bd}, the network $Q_2$ has at most
\begin{align}
	\begin{split}
	&109 (\rsum+\bsum+1)^{3+\rsum} N (m+5)
	\end{split}
	\label{eq.nr_param_Q2}
\end{align}
non-zero parameters. 

To obtain a network reconstruction of the function $f$, it remains to scale and shift the output entries. This is not entirely trivial because of the bounded parameter weights in the network. Recall that $B= \lceil 2K e^r \rceil .$ The network $x \mapsto BM_1^{r}M_2^{r_\prime}x$ is in the class $\mF(3, (1,  M_1^{r}M_2^{r_\prime} , 1, \lceil 2Ke^r\rceil,1))$ with shift vectors $\bv_j$ are all equal to zero and weight matrices $W_j$ with all entries equal to one. Because of $N \geq  (K+1)e^\rsum,$ the number of parameters of this network is bounded by $2M_1^{r}M_2^{r_\prime} +2\lceil 2Ke^r\rceil \leq 6 N$. This shows existence of a network in the class $\mF(4, (1, 2, 2M_1^{r}M_2^{r_\prime} ,2, 2\lceil 2Ke^r \rceil,  1)) $ computing $a \mapsto BM_1^{r}M_2^{r_\prime}(a-c)$ with $c:=1/(2M_1^{r}M_2^{r_\prime}).$ This network computes in the first hidden layer $(a-c)_+$ and $(c-a)_+$ and then applies the network $x \mapsto BM_1^{r}M_2^{r_\prime}x$ to both units. In the output layer, the second value is subtracted from the first one. This requires at most $6+12N$ active parameters. 

Because of \eqref{eq.Q2_estimate} and \eqref{eq.kernel_sum_one}, there exists a network $Q_3$ in
\begin{align*}
	 \mF\big((m+13) + L^*, (\rsum, 6(\rsum+\lceil \bsum\rceil)N, \ldots, 6(\rsum+\lceil \bsum\rceil)N,1)\big)
\end{align*}
such that 
\begin{flalign*}
	\Bigg | Q_3(\bu, \bv) - \sum_{ \bu_{\bell^\1}, \bv_{\bell^\2}  \in \bD(M) }  P^{\beta, \beta_\prime} f (\bu, \bv) &\Big(\prod_{j=1}^{r} (1/M_1 - |u^{(j)} - u^{(j)}_{\bell^\1}|)_+\Big)&
    \\
    &\Big(\prod_{j=1}^{r_\prime} (1/M_2 - |v^{(j)} - v^{(j)}_{\bell^\2}|)_+\Big) \Bigg|&
\end{flalign*}
\begin{flalign*}
	&\leq (2K+1) M_1^{r}M_2^{r_\prime} (1+ \mathsf{r}^2_{\mathrm{sum}} + \beta^2_{\mathrm{sum}}) (2e)^\rsum 2^{-m},\ \  \text{for all} \ (\bu,\bv) \in [0,1]^\rsum.&
\end{flalign*}
With \eqref{eq.nr_param_Q2}, the number of non-zero parameters of $Q_3$ is bounded by 
\begin{align*}
	109 (\rsum+\bsum+1)^{3+\rsum} N (m+6).
\end{align*}
Observe that  by construction $\mathsf{M} = (M_1+1)^{r} (M_2+1)^{r_\prime} \leq N \leq (3M_1)^{r}(3M_2)^{r_\prime} = 3^\rsum M^{\wtilde{r}}$ and hence $M^{-\wtilde{\beta}}\leq N^{-\wtilde{\beta}/\wtilde{r}} 3^{\rsum \wtilde{\beta}/\wtilde{r}}.$ Together with Lemma \ref{lem.Hoeld_approx}, the result follows. 
\end{proof}

\subsection{Embedding properties of neural network function classes}
\label{sec.embedding_props_NNs}
We denote $\F(L,\bm{p})$ as the class of neural networks with $L$ hidden layers and $\bm{p} \in \mathbb{N}^{L+2}$ nodes per layer. The class $\F(L,\bm{p})$ is subset of $\F(L,\bm{p})$ with the sparsity parameter $s$.

For the approximation of a function by a network, we first construct smaller networks computing simpler objects. Let $\bp = (p_0, \ldots,  p_{L+1})$ and $\bp' = (p_0',  \ldots, p_{L+1}').$ To combine networks, we make frequent use of the following rules. 

{\itshape Enlarging:} $\mF(L, \bp, s) \subseteq \mF(L, \bq, s')$ whenever $\bp \leq \bq$ componentwise and $s\leq s'.$

{\itshape Composition:} Suppose that $f \in \mF(L, \bp)$ and  $g \in \mF(L',\bp')$ with $p_{L+1} =p_0'.$ For a vector $\bv \in \mathbb{R}^{p_{L+1}}$ we define the composed network $g \circ \sigma_{\bv}(f) $ which is in the space $\mF(L+L'+1, (\bp, p_1', \ldots, p_{L'+1}')).$ In most of the cases that we consider, the output of the first network is non-negative and the shift vector $\bv$ will be taken to be zero. 

{\itshape Additional layers/depth synchronization:} To synchronize the number of hidden layers for two networks, we can add additional layers with an identity weight matrix, such that
\begin{align}
	\mF(L, \bp,s) \subset \mF(L+q, (\underbrace{p_0,\ldots,p_0}_{q\text{ times}} , \bp), s+qp_0).
	\label{eq.add_layers}
\end{align}

{\itshape Parallelization:} Suppose that $f,g$ are two networks with the same number of hidden layers and the same input dimension, that is, $f \in \mF(L, \bp)$ and $g \in \mF(L, \bp')$ with $p_0=p_0'.$ The parallelized network $(f,g)$ computes $f$ and $g$ simultaneously in a joint network in the class $\mF(L, (p_0, p_1+p_1', \ldots, p_{L+1}+p_{L+1}')).$ 

\subsection{Technical lemmas for the proof of Theorem~\ref{thm.approx_network_One_fct2}}

We use $\mF(L,{\bf{r}})$ to denote a fully connected network with $L$ deep layers and ${\bf r} \in \nat^{L+2}$ representing the nodes in each layer.

The following technical lemmas are required for the proof of \Cref{thm.approx_network_One_fct2}. \Cref{lem1}, \Cref{lem2}, and \Cref{lem3} restate Lemma A.2, Lemma A.3, and Lemma A.4 from \cite{schmidtheiber}, respectively.

\begin{lemma}\label{lem1}
\label{lem.mult}
For any positive integer $m,$ there exists a network $\Mult_m \in \mF(m+4,(2,6,6,\ldots,6,1)),$ such that $\Mult_m(x,y) \in [0,1],$
\begin{align*}
	\big|\Mult_m (x,y) - x y \big| \leq 2^{-m}, \quad \text{for all} \ x,y \in [0,1],
\end{align*}
and $\Mult_m(0,y)=\Mult_m(x,0) =0.$
\end{lemma}

\begin{lemma}\label{lem2}
\label{lem.multi_mult}
For any positive integer $m,$ there exists a network $$\Mult_m^r \in \mF((m+5) \lceil \log_2 r \rceil , (r, 6r,6r, \ldots, 6r,1))$$ such that $\Mult_m^r \in [0,1]$ and
\begin{align*}
	\Big|\Mult_m^r 
	(\bx)
	- \prod_{i=1}^r x_i \Big| \leq r^2 2^{-m}, \quad \text{for all} \ \ \bx=(x_1, \ldots,x_r)\in [0,1]^r.
\end{align*}
Moreover, $\Mult_m^r(\bx)=0$ if one of the components of $\bx$ is zero. 
\end{lemma}

The number of monomials with degree $|\balpha | < \gamma$ is denoted by $C_{r,\gamma}.$ Obviously, $C_{r,\gamma} \leq (\gamma +1)^r$ since each $\alpha_i$ has to take values in $\{0,1,\ldots,\lfloor \gamma \rfloor\}.$
\begin{lemma}\label{lem3}
\label{lem.monomials}
For $\gamma >0$ and any positive integer $m,$ there exists a network 
\begin{align*}
	\Mon_{m,\gamma}^r \in \mF \big( 1+(m+5) \lceil \log_2 (\gamma \vee 1) \rceil , (r, 6\lceil \gamma\rceil C_{r,\gamma}, \ldots, 6\lceil \gamma\rceil C_{r,\gamma}, C_{r,\gamma}) \big),
\end{align*}
such that $\Mon_{m,\gamma}^r \in [0,1]^{C_{r,\gamma}}$ and
\begin{align*}
	\Big| \Mon_{m, \gamma}^r(\bx)  -  (\bx^{\balpha})_{|\balpha| < \gamma} \Big| _\infty \leq \gamma^2 2^{-m}, \quad \text{for all} \ \bx \in [0,1]^r.
\end{align*}
\end{lemma}

\end{document}